\documentclass{amsart}
\usepackage[T1]{fontenc}
\usepackage[english]{babel}
\usepackage{amsmath,amsthm,amssymb, latexsym,geometry,stackrel,enumerate}
\usepackage[all]{xy}
\usepackage{hyperref}

\usepackage[dvips]{graphicx}
\DeclareGraphicsExtensions{.bmp,.png,.pdf,.jpgg}
\usepackage[all]{xy}
\usepackage{verbatim}
\usepackage[mathcal]{euscript}
\usepackage{epsfig}
\usepackage{multicol}
\usepackage{color}
\usepackage{verbatim}
\usepackage{graphicx,float}
\usepackage{tikz-cd}

\date{\today}

\swapnumbers
\theoremstyle{plain}
\newtheorem{teo}{Theorem}[section]
\newtheorem{lema}[teo]{Lemma}
\newtheorem{prop}[teo]{Proposition}
\newtheorem{coro}[teo]{Corollary}

\theoremstyle{definition}
\newtheorem{defi}[teo]{Definition}
\newtheorem{obs}[teo]{Remark}
\newtheorem{obss}[teo]{Remarks}
\newtheorem{ejem}[teo]{Example}

\usepackage{url}

\def\Ext{\mathop{\rm Ext}\nolimits}
\def\Tor{\mathop{\rm Tor}\nolimits}
\def\Ho{\mathop{\rm Hom}\nolimits}

\def\RMod{\mathop{\rm R\text{-} Mod}\nolimits}

\def\pdR{\mathop{\rm pd}_{R}\nolimits}
\def\id{\mathop{\rm id}\nolimits}

\def\wdR{\mathop{\rm wd}_{R}\nolimits}

\def\nid{\mathop{\rm FP_{n}\text{-}id}\nolimits}
\def\nfd{\mathop{\rm FP_{n}\text{-}fd}\nolimits}
\def\npd{\mathop{\rm FP_{n}\text{-}pd}\nolimits}
\def\glD{\mathop{\rm glD}\nolimits}

\def\fppd{\mathop{\rm FP_{n}\text{-}pd}\nolimits}

\def\1id{\mathop{\rm (1,0)\text{-}id}\nolimits}
\def\lamd{\mathop{\rm \lambda\text{-}dim}\nolimits}

\def\pd2{\mathop{\rm FP_{2}\text{-}pd}\nolimits}
\def\2gpd{\mathop{\rm FP_2\text{-}PD}\nolimits}
\def\fpD{\mathop{\rm fpD}\nolimits}
\def\igpd{\mathop{\rm FP_i\text{-}PD}\nolimits}
\def\unogpd{\mathop{\rm FP_1\text{-}PD}\nolimits}
\def\nmenos1gpd{\mathop{\rm FP_{n-1}\text{-}PD}\nolimits}
\def\nm1pd{\mathop{\rm FP_{n-1}\text{-}pd}\nolimits}
\def\fpd{\mathop{\rm fpd}\nolimits}

\def\gid{\mathop{\rm FP_n\text{-}ID}\nolimits}

\def\fgpd{\mathop{\rm f.FP_n\text{-}PD}\nolimits}
\def\gpd{\mathop{\rm FP_n\text{-}PD}\nolimits}

\newcommand{\ts}[1]{\normalfont{\textsf{#1}}}
\newcommand{\K}{\ts k}

\title[Some remarks about $FP_{n}$-projective and $FP_{n}$-injective modules]{Some remarks about $FP_{n}$-projective and $FP_{n}$-injective modules}

\author[V. Gubitosi]{Viviana Gubitosi}
\address{V. Gubitosi \newline Instituto de Matem\'{a}tica y Estad\'{\i}stica Rafael Laguardia, Facultad de Ingenier\'{\i}a - UdelaR, Montevideo, Uruguay, 11200 }
\email{gubitosi@fing.edu.uy}

\author[R. Parra]{Rafael Parra}
\address{R. Parra \newline Instituto de Matem\'{a}tica y Estad\'{\i}stica Rafael Laguardia, Facultad de Ingenier\'{\i}a - UdelaR, Montevideo, Uruguay, 11200 }
\email{rparra@fing.edu.uy}


\begin{document}
\maketitle

\begin{abstract}
In this paper, we  give some new characterizations of \( FP_n \)-projective modules and strong $n$-coherent rings. Some known results are extended and some new characterizations of the \( FP_n \)-injective global dimension in terms of \( FP_n \)-projective modules are obtained.  Using the \( FP_n \)-projective dimension of a module defined by Ouyang, Duan and Li in \cite{Ouy} we introduce a  slightly different  \( FP_n \)-projective global dimension over the ring which measures how far away the ring is from being Noetherian. This dimension agrees with the $(n,0)$-projective global dimension of \cite{Ouy} when the ring is strong $n$-coherent.

\end{abstract}

\section*{Introduction}
 The class of modules with a vanishing property with respect to the class of finitely presented modules and the functors $\Ext^1_R(-,-)$ and  $\Tor_1^R(-,-)$  have been extensively explored by many authors.  For example, Maddox \cite{Mad} and Stenstr\"om \cite{Sten} introduced \( FP \)-injective modules as $R$-modules $M$ for which \( \Ext^1_R(F, M) = 0 \) for all finitely presented \( R \)-modules \( F \), and examined them over coherent rings. In the literature, $FP$-injective modules are also known as absolutely pure modules and are often regarded as dual analogs of flat modules.

In 2005, Mao and Ding \cite{MD1} utilized the concept of \( FP \)-injective modules to define \( FP \)-projective modules. An \( R \)-module \( P \) is called \( FP \)-projective if \( \Ext_R^1(P, M) = 0 \) for any \( FP \)-injective \( R \)-module \( M \). If we denote by $\mathcal{FP}$-$Proj(R)$ and by   $\mathcal{FP}$-$Inj(R)$  the class of all \(\mathcal{FP}\)-projective $R$-modules and all  \(\mathcal{FP}\)-injective $R$-modules respectively, it is known that the pair ($\mathcal{FP}$-$Proj(R)$, $\mathcal{FP}$-$Inj(R)$), forms a complete cotorsion pair which is generated by the representative set of all finitely presented \( R \)-modules \cite{Trlifaj}. In addition, this cotorsion pair is hereditary if and only if the ring \( R \) is coherent \cite{MD2}.

Many homological results can be generalized using finitely \( n \)-presented modules instead of finitely presented modules.  In this way, \( FP \)-injective modules are replaced by \( FP_n \)-injective modules in \cite{Zhou}.  These modules are particularly useful for characterizing strong \( n \)-coherent rings, that is, rings in which every finitely \( n \)-presented module is also finitely \( (n{+}1) \)-presented. Taking $n=1$ finitely presented modules and  \( FP \)-injective modules are recovered. In the same way Mao and Ding \cite{MD4} defined the  class of \( FP_n \)-projective modules. Again, if $\mathcal{FP}_n$-$Proj(R)$ and $\mathcal{FP}_n$-$Inj(R)$ denote the classes of all \(FP_n\)-projective and \(FP_n\)-injective $R$-modules respectively, the pair  ($\mathcal{FP}_n$-$Proj(R)$, $\mathcal{FP}_n$-$Inj(R)$) forms a complete cotorsion pair that is hereditary if and only if \( R \) is a strong $n$-coherent ring \cite{BP}. Recently,   cotorsion pairs have been playing an important role in producing approximations and the existence of such approximations is a prerequisite for computing relative dimensions. 
 
Intricately connected to the concept of finitely \( n \)-presented module is found the concept of \( n \)-pure exact sequence. In \cite{MD4}, Mao and Ding characterized strong \( n \)-coherent rings studying the relationship between \( FP_n \)-injectivity and \( n \)-purity. In this work, we use the notion of $n$-pure exact sequence to demonstrate that $FP_n$-injective modules can be considered as dual analogs of $FP_n$-flat modules. This approach also allows us to extend analogous notions to those in \cite{MD1} by replacing \( FP \)-projective modules with \( FP_n \)-projective modules.

Relative homological algebra, initiated by Auslander and Buchweitz~\cite{AB}, studies homological dimensions obtained by replacing the class of projective or injective modules with certain subcategories. Since then, these ideas have been extensively developed, with many authors studying homological dimensions defined via alternative classes of modules. In this context, Ouyang, Duan, and Li~\cite{Ouy} introduced the $FP_n$-projective dimension of an $R$-module and the corresponding $FP_n$-projective global dimension of the ring~$R$. More generally, Angeleri Hügel and Mendoza~\cite{AM} studied relative homological dimensions in the setting of cotorsion pairs with applications to tilting theory and the finitistic dimension conjecture.

Throughout this paper, $R$ denotes an associative ring with a unit and unless otherwise specified all modules considered will be left $R$-modules.
This paper is organized as follows. Section 1 is devoted to recalling the concept of  \( FP_n \)-projective module and to give new characterizations. In section 2 we work with \( n \)-pure exact sequences and  \( FP_n \)-flat modules. We also introduce the class  of almost  \( FP_n \)-injective  modules and we give some characterizations of strong $n$-coherent rings. In  Section 3, we deal with the \( FP_n \)-projective dimension of a $R$-module \( M \) defined in \cite{Ouy}. Using this dimension we characterize \( n \)-von Neumann regular rings and (strong) \( n \)-coherent rings. Additionally, we introduce the corresponding global dimension over the ring \( R \), which differs from the one presented in \cite{Ouy}. We compare this dimension with other well-known dimensions, such as the global dimension and the \( \lambda \)-dimension. Motivated by classical results on projective modules over hereditary rings, we investigate when the class of \( FP_n \)-projective modules is closed under submodules and demonstrate, using the \( FP_n \)-projective dimension, that this occurs for a generalized version of hereditary rings, referred to as \( FP \)-hereditary rings. Section 4 deals with strong \( n \)-coherent and self \( FP_n \)-injective rings. With this assumption, we establish connections between the projective and \( FP_n \)-flat dimension of a  \( FP_n \)-projective module.  Section 5 is dedicated to compute the \( FP_n \)-injective global  dimension using \( FP_n \)-projective modules. As a consequence,  we characterize \( (n,d) \)-rings. In particular, characterizations of \( n \)-von Neumann regular rings and \( n \)-hereditary rings are presented. Finally, in section 6, some
applications to subprojectivity domains, the study of the $CF$-conjecture and trace modules in $FP_n$-injective envelopes are indicated.

\section{ \( FP_n \)-projective modules}

Let $n$ be a non-negative integer $(n\geq 0)$. According to \cite[Section 1]{Costa}, an $R$-module $M$ is said to be \textit{finitely $n$-presented} if there is an exact sequence:
$$F_n\rightarrow F_{n-1}\rightarrow \cdots \rightarrow F_1 \rightarrow F_0 \rightarrow M \rightarrow 0$$ where the $F_i$ are finitely generated and projective (or free) $R$-modules, for every $0\leq i \leq n$. This exact sequence is referred to as a finite $n$-presentation of $M$.

The class of all finitely $n$-presented $R$-modules is denoted by $\mathcal{FP}_n(R)$. In particular, $\mathcal{FP}_0(R)$ is the class of all finitely generated $R$-modules, and $\mathcal{FP}_1(R)$ is the class of all finitely presented $R$-modules. $\mathcal{FP}_{\infty}(R)$ represents the class of all finitely  $\infty$-presented $R$-modules, which consist of modules that have a resolution by finitely generated and projective (or free) modules. For convenience, we define $\mathcal{FP}_{-1}(R)$ to be the entire class of $R$-modules.

Additionally, we have a chain of inclusions:
$$\mathcal{FP}_0(R)\supseteq \mathcal{FP}_1(R)\supseteq \cdots \supseteq \mathcal{FP}_n(R)\supseteq \cdots\supseteq\mathcal{FP}_{\infty}(R).$$

For any class $\mathcal{C}\subseteq \rm R$-$\rm Mod $, we denote by $\mathcal{C}^{\perp} $ (or $^{\perp}\mathcal{C}$)  the right (or  left, respectively) orthogonal complement of $\mathcal{C}$; i.e.,

\begin{center}
$\mathcal{C}^{\perp}:= \{\ X\in \rm R\text{-}Mod \mid \Ext^1_R(C,X)=0 \text{ for all $C $}\in \mathcal{C}   \}$

$^{\perp}\mathcal{C}:= \{\ X\in\rm R\text{-}Mod \mid \Ext^1_R(X,C)=0 \text{ for all $C $}\in \mathcal{C}   \}$.
\end{center}

The classes $^{\perp}\mathcal{C}$ and $\mathcal{C}^{\perp}$ are closed under direct summands and extensions. Additionally, $^{\perp}\mathcal{C}$ is closed under direct sums and contains all projective modules, while $\mathcal{C}^{\perp}$ is closed under direct products and contains all injective modules.  If $\mathcal{C}$ consists of all finitely $n$-presented $R$-modules, then the class $\mathcal{C}^{\perp}$ is exactly the class of all $FP_n$-injective $R$-modules, i.e., an $R$-module $M$ is \textit{$FP_n$-injective} if $\Ext_R^1(F,M) = 0$ for all finitely $n$-presented $R$-module $F$ (this may include the case $n=\infty$). With this definition, an $R$-module $M$ is injective if and only if it is $FP_0$-injective. It is considered $FP$-injective (or absolutely pure) if and only if it is $FP_1$-injective. The class of all $FP_n$-injective $R$-modules is denoted by $\mathcal{FP}_n$-$Inj(R)$. 

\begin{defi}\cite[Definition 3.1]{MD4}.
   An $R$-module $P$ is said to be an \textit{$FP_n$-projective module} if it satisfies $\Ext^1_R(P,M)=0$ for all $FP_n$-injective $R$-modules $M$. 
\end{defi} 

The $FP_n$-projective modules generalize $FP$-projective modules, which have been studied in \cite{MD1}, and coincide with them when $n=1$. However, there are other generalizations of $FP$-projective modules in the literature. For example, Mao and Ding define $n$-$FP$-projective modules in \cite{MD2}, where $FP$-projective modules are recovered by taking $n=0$. For further information about $FP_n$-injective and $FP_n$-projective modules the reader is referred to   \cite{Zhou, MD4, Zhu, Ouy, BP, Zhu2}.

The class of all $FP_n$-projective $R$-modules is denoted by $\mathcal{FP}_n$-$Proj(R)$. According to \cite[Corollary 3.2.4]{Trlifaj} $\mathcal{FP}_n$-$Proj(R)$ consists of all direct summands of $\mathcal{S}$-filtered modules, where $\mathcal{S}$ is the representative set of all finitely $n$-presented $R$-modules. Equivalently, an $R$-module $M$ is $FP_n$-projective if and only if $M$ is a direct summand of an $R$-module $N$ such that $N$ is a union of a continuous chain $(N_{\alpha} : \alpha < \lambda)$, for some cardinal $\lambda$, with $N_0 = 0$ and $N_{\alpha +1}/N_{\alpha} \in \mathcal{FP}_n(R)$ for all $\alpha < \lambda$. \\
 
 Since every $FP_n$-injective module is also $FP_{n+1}$-injective, we obtain the following chain of inclusions:
\begin{center}
 $\mathcal{FP}_1$-$Proj(R) \supseteq \cdots \supseteq \mathcal{FP}_n$-$Proj(R)\supseteq \cdots\supseteq\mathcal{FP}_{\infty}$-$Proj(R),$
\end{center}
from where it immediately follows that $\mathcal{FP}_{\infty}$-$Proj(R) \subseteq \bigcap_{n\geq 0}\mathcal{FP}_{n}$-$Proj(R).$ 

\begin{obs}\cite[Theorem 3.9]{MD4}\label{suficientes proy} For any ring $R$, the pair of classes $(\mathcal{FP}_n$-$Proj(R), \mathcal{FP}_n$-$Inj(R))$ forms a complete cotorsion pair. Consequently, every $R$-module can be embedded into an $FP_n$-injective $R$-module in such a way that the quotient module is $FP_n$-projective, and every $R$-module is the quotient of some $FP_n$-projective $R$-module by its $FP_n$-injective submodule. Furthermore, $\mathcal{FP}_n$-$Proj(R)$ is closed under extensions, direct sums, direct summands and filtrations; while the class $\mathcal{FP}_n\text{-}\mathrm{Inj}(R)$ is closed under 
extensions, direct products, and direct summands.  If $n \geq 2$, then $\mathcal{FP}_n\text{-}\mathrm{Inj}(R)$ is also closed 
under direct limits. In particular, is closed under direct sums.
\end{obs}

Recall that a cotorsion pair \( (\mathcal{A}, \mathcal{B}) \) is \emph{hereditary} if 
\(\Ext^i_R(A,B)=0\) for all \(i \geq 2\), whenever \(A \in \mathcal{A}\) and \(B \in \mathcal{B}\).

\begin{obs}\label{par hereditario} \cite[Theorem~4.1]{MD4} and~\cite[Theorem~5.5]{BP}. A ring \( R \) is  strong left \( n \)-coherent if and only if
the cotorsion pair \( (\mathcal{FP}_n\text{-}Proj(R), \mathcal{FP}_n\text{-}Inj(R)) \) is hereditary. This is also equivalent to the fact that the class \( \mathcal{FP}_n\text{-}Inj(R) \) is  coresolving (that is, it contains the injective modules and it is closed under extensions and cokernels of monomorphisms). 

   Dually, the cotorsion pair \( (\mathcal{FP}_n\text{-}Proj(R), \mathcal{FP}_n\text{-}Inj(R)) \) is hereditary if and only if the class  \( \mathcal{FP}_n\text{-}Proj(R) \) is resolving (that is, it contains the projective modules and it is closed under extensions and kernels of epimorphisms).
\end{obs}

From now on we assume $n\geq 1$. \\
We start rewriting Theorem 2.1 of \cite{Ouy} and Theorem 2.6 of \cite{Zhu} in terms of the class of $FP_n$-projective modules.

\begin{lema} \label{FP-proj lemma}
Let $R$ be a ring and $P$ an finitely generated $R$-module. Then the following conditions are equivalent.
\begin{enumerate}
    \item $P$ is finitely $n$-presented.
    \item $P$ is $FP_n$-projective.
    \item $P$ is finitely $(n-1)$-presented $FP_n$-projective.
\end{enumerate}
\end{lema}
\qed

Observe that any projective (or free) module belongs to $ \mathcal{FP}_{n}$-$Proj(R)$, for all $n\geq 1$. The following example, which follows from the previous lemma, shows that there exist modules which are $FP_n$-projective for all $n\geq 1$ but they are not projective.

\begin{ejem} Let $A$ be a commutative ring and consider $R:=A\ltimes A$, the trivial ring extension of $A$ by $A$. Remember that the trivial ring extension of a ring $A$ by an $A$-module $E$,
also called the idealization of $E$ over $A$, is the ring $R := A \ltimes E$ whose underlying group is $A \times E$ with multiplication given by $(a_1, e_1)(a_2, e_2)=(a_1a_2, a_1e_2 +a_2e_1)$.

Consider the $R$-module $I:= 0\ltimes A$. By \cite[Proposition 2.3]{Mahdou2}, $I$ is a finitely $n$-presented $R$-module (for each positive integer $n$) which is not projective. Thus, $I$ is $FP_n$-projective for each positive integer $n$. Therefore, the inclusion $Proj(R) \varsubsetneq \bigcap_{n\geq 0}\mathcal{FP}_{n}$-$Proj(R)$ is strict.
\end{ejem}

\begin{ejem}  Let $\K$ be a field and consider the polynomial ring $R$
$$R=   \dfrac{\K[\cdots,x_3, x_2,x_1,y_1,y_2,y_3, \cdots ]}{(x_{j+1}x_{j},x_1y_{1},y_{1}y_{i})_{i,j\geq 1}}.$$
According to \cite[Example 1.4]{BP}, the ideal $(x_i)$ is in $\mathcal{FP}_i(R)$, but not in $\mathcal{FP}_{i+1}(R)$. Now, by Lemma \ref{FP-proj lemma}
we get that $(x_i)\in\mathcal{FP}_i$-$Proj(R)\setminus \mathcal{FP}_{i+1}$-$Proj(R)$ for $i\geq 1$.
\end{ejem}

Recall that a class \(\mathcal{C}\) of \(R\)-modules is closed under kernels of epimorphisms if for every short exact sequence of \(R\)-modules \(0 \rightarrow A \rightarrow B \rightarrow C \rightarrow 0\) with \(B\) and \(C\) in \(\mathcal{C}\), then \(A\) is also in \(\mathcal{C}\). In general, for any ring \(R\), the class $\mathcal{FP}_n$-$Proj(R)$ is not closed under kernels of epimorphisms. Remark~\ref{par hereditario} specifies when this property holds true.

By \cite[Lemma 3.3]{MD4}, if  \(0 \rightarrow A \rightarrow B \rightarrow C \rightarrow 0\) is a short exact sequence of \(R\)-modules such that \(C\) is \( FP_{n+1} \)-projective and \(B\) is \( FP_n \)-projective, then \(A\) is \( FP_n \)-projective. We use this fact in the following proposition. 

\begin{prop}
Let $R$ be a ring and let $M$ and $P$ be $FP_n$-projective $R$-modules such that $M + P$ is also $FP_n$-projective. Then $M \cap P$ is $FP_{n-1}$-projective.
\end{prop}
\begin{proof}
Since $M \oplus P$ is $FP_n$-projective, by Lemma \ref{FP-proj lemma}, it is also $FP_{n-1}$-projective. Now use \cite[Lemma 3.3]{MD4} and the exact sequence $0 \rightarrow M \cap P\ \rightarrow M \oplus P \rightarrow M + P \rightarrow 0$.
\end{proof}

Recall that an $R$-module $P$ is said to be projective with respect to a short exact sequence $\mathcal{S}$ of $R$-modules if $\Ho_R(P, \mathcal{S})$ is exact.

\begin{prop}\label{FP-proj prop} Let $R$ be a ring and $P$ be an $R$-module. The following conditions are equivalent.
\begin{enumerate}
\item $P$ is $FP_n$-projective.
\item $P$ is projective with respect to every short exact sequence $0\rightarrow A \rightarrow B \rightarrow C \rightarrow 0$ of $R$-modules with  $A \in \mathcal{FP}_n$-$Inj(R)$.
\item Every short exact sequence $0\rightarrow A \rightarrow B \rightarrow P \rightarrow 0$ of $R$-modules, with $A \in \mathcal{FP}_n$-$Inj(R)$, splits.
\end{enumerate}
\end{prop}

\begin{proof}
$(1)\Rightarrow(2)$ Assume that $P$ is $FP_n$-projective. Consider a short exact sequence of $R$-modules  $0\rightarrow A \rightarrow B \rightarrow C \rightarrow 0$, where $A$ is $FP_n$-injective. By definition, we have $\Ext^1_R(P,A)= 0$. From the induced exact sequence
$$\cdots\rightarrow \Ho_R(P,B)\rightarrow \Ho_R(P,C)\rightarrow \Ext^1_R(P,A)\rightarrow \Ext^1_R(P,B)\rightarrow \cdots $$
we conclude that  $\Ho_R(P,B)\rightarrow \Ho_R(P,C)\rightarrow 0$ is exact. Hence, $P$ is projective with respect to the short exact sequence  $0\rightarrow A \rightarrow B \rightarrow C \rightarrow 0$.\\

$(2)\Rightarrow(3)$ and  $(3) \Rightarrow (1)$ are clear. 
\end{proof}

\section{$n$-pure exact sequences and $n$-pure modules} 

A short exact sequence $\mathcal{S}$ of $R$-modules $0\rightarrow F \rightarrow E \rightarrow G \rightarrow 0$ is called \emph{pure} if and only if $\Ho_R(M, \mathcal{S})$ is exact for each finitely $1$-presented  $R$-module $M$. When $E$ is flat, it is well known that $G$ is flat if and only if $\mathcal{S}$ is pure \cite[36.6]{Wis}. A submodule $A$ of an $R$-module $B$ is said to be a \textit{pure submodule} if the induced map $\Ho_R(M, B)\rightarrow \Ho_R(M, B/A)$ is surjective for all finitely $1$-presented $R$-module $M$. Recall that an $R$-module $M$ is $FP_1$-injective if and only if it is a pure submodule of every overmodule, (i.e., every module containing it as a submodule).

Short exact sequences arising from the canonical presentation of a direct limit form an important class of examples of short pure exact sequences. That is, let $(M_i, f_{ij})_{i,j\in I}$ be a direct system of modules and consider its direct limit $\lim_{\rightarrow I} M_i$. The canonical presentation
$$ 0 \rightarrow \rm Ker(\pi) \rightarrow \bigoplus_{i\in I}M_i \xrightarrow{\pi}  \lim_{\rightarrow I}M_i \rightarrow 0$$
of $\lim_{\rightarrow I}M_i$ is an example of a short pure exact sequence \cite{Trlifaj}.\\

As defined in \cite{MD4} a short exact sequence $\mathcal{S}$ of $R$-modules is called {\emph{$n$-pure}} if every finitely $n$-presented  $R$-module is projective with respect to this sequence $\mathcal{S}$. Analogously,  a submodule $A$ of an $R$-module $B$ is said to be a \textit{n-pure submodule} (or $A$ is $n$-pure in $B$ for short) if the induced map $\Ho_R(M, B)\rightarrow \Ho_R(M, B/A)$ is surjective for all finitely $n$-presented $R$-module $M$.\\

A very strong and useful  result about $n$-pure exact sequences was obtained recently by Tan, Wang and Zhao. They stated that an exact sequence of $R$ modules $0\rightarrow A \rightarrow B \rightarrow C \rightarrow 0$ is  $n$-pure if and only if the induced sequence $0 \rightarrow M \otimes_R A \rightarrow M \otimes_R B \rightarrow M \otimes_R C \rightarrow 0$ is exact for every finitely $n$-presented right $R$-module $M$. See \cite[Theorem 2.5]{TWZ}.\\

This result allows us to easily generalize a well known property of pure submodules. \\

\begin{prop} Let \( A_1, A_2 \) be submodules of an \( R \)-module \( A \) and $n \geq 1$ an integer. 
\begin{enumerate}
    \item  If $A_1 \subset A_2$ and $A_2$ is $n$-pure in $A$, then $A_2/A_1$ is $n$-pure in $A/A_1$. 
    \item   If $A_1 \subset A_2$ and $A_1$ is $n$-pure in $A$ and $A_2/A_1$ is $n$-pure in $A/A_1$ then $A_2$ is $n$-pure in $A$.
    \item  If \( A_1 + A_2 \) and \( A_1 \cap A_2 \) are \( n \)-pure in \( A \), then \( A_1 \) and \( A_2 \) are \( n \)-pure in \( A \).
\end{enumerate}
\end{prop}

\begin{proof}
For $n=1$ it is known. See \cite[Ex. 4.30-4.31]{Lam2}. For $n>1$ observe that it follows in the same way as the case \( n = 1 \), using \cite[Theorem~2.5]{TWZ}.
\end{proof}

Recall that an $R$-module $M$ is called \textit{$FP_n$-flat} if $\Tor^R_1(F,M)=0$ for all finitely $n$-presented right $R$-module $F$ (this may include the case $n=\infty$).  Note that the class of $FP_0$-flat modules coincide with the class of flat modules.  The class of all $FP_n$-flat $R$-modules is denoted by $\mathcal{FP}_n$-$Flat(R)$.\\

Another consequence of \cite[Theorem 2.5]{TWZ} is the following lemma.

\begin{lema}\label{pur2}  Let $R$ be a ring. An $R$-module $C$ is $FP_n$-flat if and only if every exact sequence $0 \rightarrow A \rightarrow B \rightarrow C \rightarrow 0$ of $R$-modules is $n$-pure.
\end{lema}

\begin{proof}
Assume that $C$ is $FP_n$-flat. Then, for any finitely $n$-presented right $R$-module $M$, we have the following exact sequence: 
\[ \Tor^R_1(M,C) \rightarrow M \otimes_R A \rightarrow M \otimes_R B \rightarrow M \otimes_R C \rightarrow 0. \]
Since $C$ is $FP_n$-flat, $\Tor^R_1(M,C) = 0$, and therefore the sequence $0 \rightarrow A \rightarrow B \rightarrow C \rightarrow 0$  is $n$-pure by \cite[Theorem 2.5]{TWZ}.

Conversely,  choose an exact sequence $0 \rightarrow A \rightarrow B \rightarrow C \rightarrow 0$ of $R$-modules with $B$  free. For every finitely $n$-presented right $R$-module $M$,  we have the following exact sequence: 
\[ 0 = \Tor^R_1(M,B) \rightarrow \Tor^R_1(M,C) \rightarrow M \otimes_R A \rightarrow M \otimes_R B \rightarrow M \otimes_R C \rightarrow 0. \]
Since the sequence \[ 0 \rightarrow M \otimes_R A \rightarrow M \otimes_R B \rightarrow M \otimes_R C \rightarrow 0 \] is exact by \cite[Theorem 2.5]{TWZ},  it follows that $\Tor^R_1(M,C) = 0$. Therefore,  $C$ is $FP_n$-flat.
\end{proof}

Additionally, since every \( R \)-module can be expressed as the quotient of a free module, we have:

\begin{coro}
    The following conditions are equivalent.
    \begin{enumerate}
        \item \( C \) is \( FP_n \)-flat.
        \item Every exact sequence \( 0 \rightarrow A \rightarrow B \rightarrow C \rightarrow 0 \) of \( R \)-modules is \( n \)-pure.
        \item There exists an \( n \)-pure exact sequence \( 0 \rightarrow A \rightarrow B \rightarrow C \rightarrow 0 \) of \( R \)-modules where \( B \) is \( FP_n \)-flat.
    \end{enumerate}
\end{coro}
\qed

It follows directly from Lemma \ref{pur2} the following Corollary.

\begin{coro}\label{pur3}
    Let \( R \) be a ring and consider an exact sequence \(\mathcal{S}: 0 \rightarrow A \rightarrow B \rightarrow C \rightarrow 0 \)  of \( R \)-modules, where $B$ is \( FP_n \)-flat. Then, $\mathcal{S}$ is $n$-pure if and only if $C$ is \( FP_n \)-flat.
\end{coro}
\qed

It is well known that a ring \( R \) is von Neumann regular if and only if every finitely presented \( R \)-module is projective, or equivalently, if every short exact sequence of \( R \)-modules is pure. A ring \( R \) is said to be left \emph{\( n \)-von Neumann regular} if every finitely \( n \)-presented \( R \)-module is projective. According to \cite[Theorem 3.9]{Zhu}, this is equivalent to the condition that every right \( R \)-module is \( FP_n \)-flat. As an immediate consequence of the previous lemma, we obtain:

\begin{coro}\label{n-von Neumann y sec n-puras: el (2.5)}

    A ring \( R \) is left \( n \)-von Neumann regular if and only if every short exact sequence of (left or right) \( R \)-modules is \( n \)-pure. 
\end{coro}
\qed

\begin{prop}
    Let \( R \) be a ring and \( M_1, M_2 \) be two $n$-pure submodules of an \( R \)-module \( M \). If  \( M_1 + M_2 \) is \( FP_n \)-flat then,  \( M_1 \cap M_2 \) is \( n \)-pure in \( M \).
\end{prop}

\begin{proof}
    Consider the exact commutative diagram:
    \[
    \xymatrix{
    0 \ar[r] & M_1 \cap M_2 \ar[d] \ar[r]^\alpha & M_1 \oplus M_2 \ar[d] \ar[r]^\beta & M_1 + M_2 \ar[d] \ar[r] & 0 \\
    0 \ar[r] & M \ar[r]^\gamma & M \oplus M \ar[r]^\delta & M \ar[r] & 0
    }
    \]
    where \(\alpha(x) = (x, -x)\), \(\beta(x, y) = x + y\), \(\gamma(m) = (m, -m)\) and \(\delta(m, m_0) = m + m_0\). Applying the functor $X\otimes_R-$ to the above  diagram with \( X \) any finitely \( n \)-presented right \( R \)-module  and using Lemma \ref{pur2}, we obtain the following exact commutative diagram:
    \[
    \xymatrix{
    0 \ar[r] &  X \otimes (M_1 \cap M_2)  \ar[d]^f \ar[r] & X \otimes(M_1 \oplus M_2)  \ar[d]^g \ar[r] & X \otimes (M_1 + M_2)  \ar[d]^h \ar[r] & 0 \\
    0 \ar[r] & X \otimes M  \ar[r] & X \otimes (M \oplus M)  \ar[r] & X \otimes M  \ar[r] & 0
    }
    \]
    The \( n \)-purity of \( M_1 \) and \( M_2 \) in \( M \) implies that the map \( g \) is injective. Therefore, \( f \) is also injective  which implies that \( M_1 \cap M_2 \) is \( n \)-pure in \( M \).
\end{proof}

\begin{lema}\label{pur}
Let \( R \) be a ring and consider an \( n \)-pure exact sequence \(\mathcal{S}:  0 \rightarrow A \rightarrow B \rightarrow C \rightarrow 0 \) of \( R \)-modules. If \( B \) is \( FP_n \)-injective then,  \( A \) is \( FP_n \)-injective. 
\end{lema}

\begin{proof}
The result follows by the same argument as in the well-known case \( n = 1 \).
\end{proof}

A direct consequence of Lemma \ref{pur} is the following corollary.

\begin{coro}\label{submodulo puro de Fpn iny es FPn iny}
 Every $n$-pure submodule of a $FP_n$-injective $R$-module is $FP_n$-injective. 
\end{coro}
\qed

\begin{teo}\label{FP-proj pureza} Let $R$ be a ring and $P$ an $R$-module. Then the following conditions are equivalent.
\begin{enumerate}
\item $P$ is $FP_n$-projective.
\item $P$ is projective with respect to every short $n$-pure exact sequence $0\rightarrow A \rightarrow B \rightarrow C \rightarrow 0$ of $R$-modules such that $B$ is $FP_n$-injective.
\end{enumerate}

 If  $\mathcal{FP}_n\text{-Proj}(R) \subseteq \mathcal{FP}_n\text{-Inj}(R)$,  these conditions are also equivalent to:
\begin{enumerate}\setcounter{enumi}{2}

\item Every short $n$-pure exact sequence $0\rightarrow A \rightarrow B \rightarrow P \rightarrow 0$ of $R$-modules  with $B\in \mathcal{FP}_n$-$Inj(R)$, splits.
\end{enumerate}
\end{teo}

\begin{proof}
$(1)\Rightarrow(2)$ Follows from Lemma \ref{pur} and Proposition \ref{FP-proj prop}.

$(2)\Rightarrow(1)$ Let $N$ be an  $FP_n$-injective $R$-module. By \cite[Theorem 2.2]{Zhu} there exists a short $n$-pure exact sequence of $R$-modules $$0\rightarrow N \rightarrow E \rightarrow K \rightarrow 0$$ where $E$ is injective.

To show that $P$ is $FP_n$-projective, we need to prove that $\Ext^1_R(P,N)= 0$.

We consider the following induced exact sequence:
$$\cdots\rightarrow \Ho_R(P,E)\rightarrow \Ho_R(P,K)\rightarrow \Ext^1_R(P,N)\rightarrow \Ext^1_R(P,E)\rightarrow \cdots $$
By assumption, $\Ho_R(P,E)\rightarrow \Ho_R(P,K)\rightarrow 0$ is exact. Also, since $E$ in injective, we have $\Ext^1_R(P,E)= 0$. Therefore, from the induced sequence, it follows that $\Ext^1_R(P,N)= 0$. Hence, $P$ is $FP_n$-projective.\\
$(2)\Rightarrow(3)$  It is clear.\\
$(3)\Rightarrow (1)$: Assume that $\mathcal{FP}_n\text{-Proj}(R) \subseteq \mathcal{FP}_n\text{-Inj}(R)$. By Remark~\ref{suficientes proy}, for any $R$-module $P$ there exists a short exact sequence $
0 \longrightarrow A \longrightarrow B \longrightarrow P \longrightarrow 0$
where $B$ is $FP_n$-projective and $A$ is $FP_n$-injective. Since $B$ is $FP_n$-injective by assumption, the sequence splits, showing that $P$ is a direct summand of $B$ and hence $FP_n$-projective.
\end{proof}

Lemma \ref{pur2} and  Corollary \ref{pur3} imply the following. 

\begin{coro} Let $R$ be a ring. If $P$ is $FP_n$-projective, then $P$ is projective with respect to every short exact sequence $0\rightarrow A \rightarrow B \rightarrow C \rightarrow 0$ of $R$-modules such that $B$ is $FP_n$-injective and $C$ is $FP_n$-flat. Moreover, if  $\mathcal{FP}_n$-$Inj(R)\subseteq \mathcal{FP}_n$-$Flat(R)$ the converse holds. 
\end{coro}

From the previous theorem, we obtain the following result.

\begin{coro} Let \( R \) be a ring. Consider the following diagram
  \[
    \xymatrix{
  & Q \ar[d]^\alpha&  & P \ar[d]^\gamma  \\
 0\ar[r] & A \ar[d] \ar[r]^f & B \ar[r]^g & C \ar[d] \ar[r] & 0 \\
   & 0  &   & 0 &  
    }
    \]
    where the row is $n$-pure exact, \( P \) is an \( FP_n \)-projective \( R \)-module and \( B \) is an \( FP_n \)-injective \( R \)-module. Then we have the following commutative diagram:
    \[
    \xymatrix{
       & 0 \ar[d] & 0 \ar[d] & 0 \ar[d] &  \\
   0 \ar[r] & A_1 \ar[d] \ar[r]  & B_1 \ar[d] \ar[r] & C_1 \ar[d] \ar[r] & 0 \\
 0\ar[r]  & Q \ar[d]^\alpha \ar[r] & Q\oplus P \ar[d]^\beta \ar[r] & P\ar@{-->}[dl]_h \ar[d]^\gamma \ar[r] & 0 \\
 0\ar[r] & A \ar[d] \ar[r]^f & B\ar[d] \ar[r]^g & C\ar[d] \ar[r] & 0 \\
   & 0  & 0 & 0 &  
    }\]
where all rows and columns are exact.
\end{coro}
\qed

\begin{obs}\label{envelopin} Recall that a pair  $(\mathcal{M}, \mathcal{C})$ in $\text{R}$-$\text{Mod}$ is a \emph{duality pair}  if it satisfies the following conditions:
\begin{enumerate}
\item $M\in\mathcal{M}$ if and only if $M^{+}=\Ho_{\mathbb{Z}}(M, \mathbb{Q}/\mathbb{Z}) \in\mathcal{C}$.
\item The class $\mathcal{C}$ is closed under direct summands and under finite direct sums.\\
\end{enumerate}

It is known that the class $\mathcal{FP}_n$-$Proj(R)$  is closed under pure quotients if and only if the pair $(\mathcal{FP}_n$-$Proj(R), \mathcal{FP}_n$-$Flat(R)^\perp)$ is a duality pair \cite[Proposition 3.7]{BEI}. Moreover, if these conditions hold,  the class  $\mathcal{FP}_n$-$Inj(R)$ is enveloping. Consequently, if $R$ is a ring where the class of $FP_n$-projective modules is closed under pure quotients, we have that $P\in\mathcal{FP}_n$-$Proj(R)$ if and only if its Pontryagin dual $P^{+}\in \mathcal{FP}_n$-$Flat(R)^\perp$.
\end{obs}

Since every pure submodule is also $n$-pure for any $n>1$, if the class of $\mathcal{FP}_n$-projective $R$-modules is closed under $n$-pure quotients, it is also closed under pure quotients. Therefore, we can characterize $\mathcal{FP}_n$-projective modules in this context using \cite[Proposition 3.7]{BEI}.

\begin{prop}
Let $n \geq 2$ and \( R \)  a ring such that the class of \( FP_n \)-projective modules is closed under $n$-pure quotients. Then, the pair $(\mathcal{FP}_n$-$Proj(R), \mathcal{FP}_n$-$Flat(R)^\perp)$ is a duality pair.
\end{prop}
\qed

We obtain the following result as an immediate consequence of the previous proposition and \cite[Theorem 2.5]{Zhu2}.

\begin{coro}
Let $n \geq 2$ and \( R \)  a ring where the class of $FP_n$-projective modules is closed under $n$-pure quotients. Then, for any $R$-module $P$, the following statements are equivalent.
\begin{enumerate}
    \item  $P$  is \( FP_n \)-projective.
    \item $ P^+$ is injective with respect to every exact sequence $0 \rightarrow A \rightarrow B \rightarrow C \rightarrow 0$ of right $R$-modules where $C$ is $FP_n$-flat.
    \item For any $FP_n$-flat right $R$-module $F$, $F$ is projective with respect to every exact sequence $0 \rightarrow  P^+ \rightarrow B \rightarrow C \rightarrow 0$ of right $R$-modules
\end{enumerate}
Moreover, if the injective envelope $E(P^+)$ of $P^+$ is $FP_n$-flat, then the above conditions are also equivalent to:
\begin{enumerate}
    \setcounter{enumi}{3}
    \item If the sequence $0 \rightarrow P^+ \rightarrow F \rightarrow L \rightarrow 0$ is exact,  where $F$ is $FP_n$-flat, then $F \rightarrow L \rightarrow 0$ is an $FP_n$-flat precover of $L$.
    \item $P^+$ is a kernel of an $FP_n$-flat precover $E \rightarrow L$ with $E$ injective.
\end{enumerate}
\end{coro}
\qed

Recall that an $R$-module $M$ is \textit{$n$-pure projective} if for any $n$-pure exact sequence
\[ 0 \rightarrow A \rightarrow B \rightarrow C \rightarrow 0 \]
of $R$-modules, the induced sequence
\[ 0 \rightarrow \Ho_R(M, A) \rightarrow \Ho_R(M, B) \rightarrow \Ho_R(M, C) \rightarrow 0 \]
is exact. See \cite[Definition 2.2]{TWZ}. Using \cite[Theorem 33.6]{Wis}, we obtain the following characterization of $n$-pure projective modules:

\begin{prop}\label{prop5}
Let $M$ be an $R$-module. Then the following statements are equivalent.
\begin{enumerate}
    \item $M$ is $n$-pure projective.
    \item Every $n$-pure exact sequence $0 \rightarrow K \rightarrow P \rightarrow M \rightarrow 0$ of $R$-modules splits.
    \item $M$ is a summand of a direct sum of finitely $n$-presented $R$-modules.
\end{enumerate}
\end{prop}
\qed

\noindent Obviously, every \( n \)-pure projective module is \( FP_n \)-projective.\\

Recall from \cite{DKM} that $R$ is called left {\emph{$n$-coherent}} if each finitely $(n-1)$-presented ideal of $R$ is finitely $n$-presented; and that $R$ is {\emph{strong left $n$-coherent}}  if each finitely $n$-presented $R$-module is finitely $(n+1)$-presented. It is known that every strong left $n$-coherent ring is left $n$-coherent. The converse is true  for $n=1$, but remains as an open question for $n\geq 2$. The 1-coherent rings are just known as coherent rings. Inspired by the work of Moradzadeh-Dehkordi and  Shojaee \cite[Theorem 3.7]{Shojaee} we show the following theorem.

\begin{teo}
    The following conditions are equivalent for a ring $R$.
    \begin{enumerate}
    \item $R$ is strong left $n$-coherent.
    \item Every finitely $(n-1)$-presented submodule of a projective (free) $R$-module is $n$-pure projective.
    \item Every finitely $(n-1)$-presented submodule of a projective (free) $R$-module is $FP_n$-projective.
    \item Every finitely $n$-presented $R$-module is an $FP_{n+1}$-projective $R$-module.
    \end{enumerate}
\end{teo}

\begin{proof}
    $(1) \Leftrightarrow (3) \Leftrightarrow (4)$ Follow from Lemma \ref{FP-proj lemma} and \cite[Theorem 2.1]{Zhu}.
    
    $(2) \Rightarrow (3)$ It is clear.
    
    $(3) \Rightarrow (2)$ Follows from Lemma \ref{FP-proj lemma} and Proposition \ref{prop5}.
\end{proof}

\noindent Similarly, with analogous arguments, we characterize $n$-coherent rings as follows:

\begin{teo}
    The following conditions are equivalent for a ring $R$.
    \begin{enumerate}
    \item $R$ is left $n$-coherent.
    \item Every finitely $(n-1)$-presented ideal of $R$ is $n$-pure projective.
    \item Every finitely $(n-1)$-presented ideal of $R$ is $FP_n$-projective.
    \end{enumerate}
\end{teo}
\qed

\begin{prop} 
        Let \( R \) be a strong right \( n \)-coherent ring, $
    0 \rightarrow M_k \rightarrow \cdots \rightarrow M_1 \rightarrow 0 $
    an exact sequence of \( FP_n \)-flat \( R \)-modules and \( N \) any finitely \( n \)-presented right \( R \)-module. Then the  sequence  $0 \rightarrow N \otimes M_k \rightarrow \cdots \rightarrow N \otimes M_1 \rightarrow 0$ is exact.
   
\end{prop}

\begin{proof}
    Set \( K := \ker(M_2 \rightarrow M_1) \). Then, we have the following two exact sequences:
    \[
    0 \rightarrow M_k \rightarrow \cdots \rightarrow M_3 \rightarrow K \rightarrow 0  \text{ \hspace*{.5cm} and  \hspace*{.5cm}}   0 \rightarrow K \rightarrow M_2 \rightarrow M_1 \rightarrow 0.
    \]

    Since \( M_1 \) is \( FP_n \)-flat, the sequence 
    $
    0 \rightarrow K \rightarrow M_2 \rightarrow M_1 \rightarrow 0 
    $
    is \( n \)-pure. Consequently, for any finitely \( n \)-presented right \( R \)-module \( N \), the sequence    $0 \rightarrow N \otimes K \rightarrow N \otimes M_2 \rightarrow N \otimes M_1 \rightarrow 0$
    is exact by \cite[Theorem 2.5]{TWZ}.

    In addition, since \( M_1 \) and \( M_2 \) are \( FP_n \)-flat, \( K \) is \(FP_n \)-flat by \cite[Corollary 2.20]{Zhu}. Thus, by induction on \( k \), the sequence  $0 \rightarrow N \otimes M_k \rightarrow \cdots \rightarrow N \otimes M_3 \rightarrow N \otimes K \rightarrow 0$ is exact. Therefore, the sequence  $0 \rightarrow N \otimes M_k \rightarrow \cdots \rightarrow N \otimes M_1 \rightarrow 0$
    is exact.
\end{proof}

\subsection{Almost $FP_n$-injective and $FP_n$-injective modules.}

The class of $FP_n$-injective modules has been extensively studied. See for example \cite{BP,Zhu}.  In this section we give some characterization of self $FP_n$-injective rings and strong left $n$-coherent rings using the class of $FP_n$-injective modules, for $n \geq 2$.  \\

Recall that $R$ is a left self $FP_n$-injective ring if $R$ is $FP_n$-injective as an $R$-module.

 Note that, by \cite[Theorem 2.2]{Zhu}, $R$ is a left self $FP_n$-injective ring if and only if it is an $n$-pure submodule of its injective envelope.  In particular, every left \( n \)-von Neumann regular ring is left self-\( FP_n \)-injective by Corollary \ref{n-von Neumann y sec n-puras: el (2.5)}.

\begin{prop}\label{autoninj FPn}  A ring $R$ is left self $FP_n$-injective if and only if for any $FP_n$-flat $R$-module $F$, there exists a short exact sequence $0 \rightarrow D \rightarrow E \rightarrow F \rightarrow 0$ of $R$-modules  where $E$ is an $FP_n$-injective module and $D$ is an $n$-pure submodule of $E$.
\end{prop}

\begin{proof}Assume that $R$ is a left self $FP_n$-injective ring and let $F$ be an $FP_n$-flat $R$-module. There exists an exact sequence of $R$-modules $$0 \rightarrow D \rightarrow E \rightarrow F \rightarrow 0$$ where $E$ is free. By Lemma \ref{pur2}, this sequence is $n$-pure. From Remark \ref{suficientes proy}, it follows that $E$ is $FP_n$-injective.

On the other hand, we obtain the following split exact sequence:  $$0\rightarrow D \rightarrow E \rightarrow R \rightarrow 0,$$  where $E$ is an $FP_n$-injective module. According to Remark \ref{suficientes proy}, this implies that R is left self $FP_n$-injective.
\end{proof}

\begin{obs}\label{almos} According to \cite{Couchot} an  $R$-module $M$ is said to be  \textit{almost $FP_1$-injective} if there exists an $FP_1$-injective $R$-module $E$ and a pure submodule $D$ such that $M$ is isomorphic to $E/D$.
For $n=1$,  Proposition \ref{autoninj FPn} states that a ring $R$ is left self $FP_1$-injective if and only if any flat $R$-module $F$  is almost $FP_1$-injective. Therefore,  we recover \cite[Proposition 3]{Couchot}.
\end{obs}

Motivated by this, we present the following definition.

\begin{defi}
An $R$-module $M$ is said to be an \textit{almost $FP_n$-injective} module if there exists an $FP_n$-injective $R$-module $E$ and an $n$-pure submodule $D$ such that $M$ is isomorphic to  $E/D$.
\end{defi}

We conclude with a new characterization of strong left $n$-coherent rings.

\begin{teo}
 A ring $R$ is strong left $n$-coherent if and only if each almost $FP_n$-injective $R$-module $M$ is $FP_n$-injective.
\end{teo}
\begin{proof}

Assume that $M$ is an almost $FP_n$-injective $R$-module over a strong left $n$-coherent ring $R$. Then, $M$ is isomorphic to  $E/D$ with $D$ an $n$-pure submodule of an $FP_n$-injective module $E$ and therefore an $FP_n$-injective module itself (by Corollary \ref{submodulo puro de Fpn iny es FPn iny}). Finally, the claim follows by observing that the class $\mathcal{FP}_n$-$Inj(R)$ is coresolving (by Remark~\ref{par hereditario}). 

For the converse, we must prove that the class $\mathcal{FP}_n$-$Inj(R)$ is coresolving. Then, it is enough to show that $\mathcal{FP}_n$-$Inj(R)$ is closed under cokernels of monomorphisms. Clearly, every such cokernel is an almost $FP_n$-injective module.

\end{proof}

It is well-known that in a coherent and self \(FP_1\)-injective ring \(R\), every flat module is \(FP_1\)-injective \cite{Sten}. The following proposition generalizes this result.

\begin{prop}
If \(R\) is a strong left \(n\)-coherent and left self \(FP_n\)-injective ring, then every \(FP_n\)-flat \(R\)-module is \(FP_n\)-injective.
\end{prop}

\begin{proof}
Let \(L\) be an \(FP_n\)-flat \(R\)-module. Consider the exact sequence \(0 \rightarrow K \rightarrow F \rightarrow L \rightarrow 0\), where \(F\) is a free \(R\)-module. Clearly, \(F\) is \(FP_n\)-injective and since \(L\) is \(FP_n\)-flat, the sequence is $n$-pure by Lemma \ref{pur2}. Consequently, \(K\) is \(FP_n\)-injective by Corollary \ref{submodulo puro de Fpn iny es FPn iny}. Finally, since the class $\mathcal{FP}_n$-$Inj(R)$ is coresolving, $L$ is  \(FP_n\)-injective.

\end{proof}

\section{$FP_n$-Projective dimension }

In \cite{MD1}, Mao and Ding defined the $FP$-projective dimension for modules and rings, denoted by $\fpd$ and $\fpD$ respectively. Their definitions measure how far away a finitely generated module is from being finitely presented and how far away a ring is from being Noetherian. With the additional assumption of coherence, they showed that
the $FP$-projective dimension has the properties that we expect of a dimension. In particular, they showed that if in addition $R$ is a left coherent ring then, $\fpD(R)= \sup \{\fpd(M)\mid M \text{ is  a left $R$-module}\}.$\

Later, Ouyang, Duan and Li generalized these dimensions in \cite{Ouy}. 

 \begin{defi}\cite[Definition 2.1]{Ouy}. The \textit{left \( FP_n \)-projective dimension} of an $R$-module \( M \), denoted by \( \npd_R(M) \), is the smallest non-negative integer \( k \) such that \( \Ext_R^{k+1}(M,N) = 0 \) for all \( FP_n \)-injective \( R \)-modules \( N \). If no such \( k \) exists, \( \npd_R(M) \) is defined to be infinite.\end{defi}

Clearly, $M$ is an $FP_n$-projective $R$-module if and only if $\npd_R(M)= 0$. In addition,  for every $R$-module $M$, $\npd_R(M)\leq \pdR(M)$ and $\nm1pd_R(M)\leq \npd_R(M)$.  The $n$-von Neumann regular rings can be characterized using the $FP_n$-projective dimension and  $FP_n$-projective modules.

\begin{prop}\label{von Neumann} The following are equivalent for a ring $R$.
\begin{enumerate}
\item $R$ is a left $n$-von Neumann regular ring.
\item $\npd_R(M)=\pdR(M)$, for any $R$-module $M$.
\item Every $FP_n$-projective $R$-module is flat.
\end{enumerate}
\end{prop}

\begin{proof}
$(1)\Rightarrow (2)$ Follows directly from \cite[Theorem 4.8]{MD4} which states that $R$ is a left $n$-von Neumann regular if and only if every $R$- module is $FP_n$-injective and therefore $\pdR(M)\leq \npd_R(M)$.\\
$(2)\Rightarrow (3)$ It is clear. \\
$(3)\Rightarrow (1)$ By Lemma \ref{FP-proj lemma}, every finitely $n$-presented module is $FP_n$-projective and by assumption, flat. Thus, every finitely $n$-presented module is projective.
\end{proof}

As usual, given a dimension over the $R$-modules we can consider the corresponding global dimension over the ring $R$.

\begin{defi}
The left \textit{$FP_n$-projective global dimension} of a ring $R$, denoted by $\gpd(R)$, is defined as follows:
$$\gpd(R)= \sup \{\fppd_R(M)\mid M \text{ is a $R$-module} \}.$$
\end{defi}

Note that in \cite{Ouy} the authors define a slightly different $FP_n$-projective global dimension of a ring $R$, called $(n,0)$-projective dimension of $R$,  by taking the above supreme over the finitely generated $R$-modules. Over strong left $n$-coherent rings both definitions agree. (See \cite[Theorem 3.1]{Ouy}). \\

Observe that $ \nmenos1gpd(R) \leq \gpd(R)$ for all $n\geq 2$.

It follows directly from the definition that rings \( R \) with \(\gpd(R) = 0\) are Noetherian rings. We will show that rings \( R \) with \(\gpd(R) = 1\) are precisely the \( FP \)-hereditary rings (see \cite{Shojaee}) that are not Noetherian.\\

 Recall that a ring $R$ is called left \textit{hereditary} if every ideal of $R$ is projective. A ring $R$ is said to be left \textit{$FP$-hereditary} if every ideal of $R$ is $FP_1$-projective \cite{Shojaee}. It is clear that  every left hereditary ring is a left $FP$-hereditary ring. In  \cite[Corollary 3.8]{Shojaee}, it was shown that every left $FP$-hereditary ring is also a left coherent ring.
 
 Non trivial examples of left $FP$-hereditary rings are the left  coherent rings where every left ideal is countably generated. See \cite[Proposition 2.3]{Posi} and \cite[Theorem 3.16]{Shojaee}. This implies that the valuation rings with only zero as  zero divisors and with a countable spectrum are also examples of $FP$-hereditary rings. See \cite[Theorem II.11]{Couchot3} and \cite[Corollary 36]{Couchot2}.\\

 We can characterize the $FP$-hereditary rings using $FP_n$-projective modules. 

\begin{prop}
Let $R$ be a ring and $n \geq 2$. The following are equivalent.
\begin{enumerate}
    \item $R$ is left $FP$-hereditary.
    \item Every ideal of $R$ is $FP_n$-projective.
    \item $\mathcal{FP}_n$-$Proj(R)$ is closed under submodules.
    \item Every $FP_n$-injective $R$-module has injective dimension at most 1.
\end{enumerate}
\end{prop}

\begin{proof}
$(1) \Rightarrow (2)$ If $R$ is a left $FP$-hereditary ring, then $R$ is a left coherent ring and \cite[Theorem 5.5]{BP} implies that $\mathcal{FP}_1$-$Inj(R) = \mathcal{FP}_n$-$Inj(R)$ for all $n > 1$. This means that every ideal in $R$ is $FP_n$-projective.

$(2) \Rightarrow (3)$ 
Assume that every ideal of $R$ is $FP_n$-projective. Thus, $R$ is $FP$-hereditary and consequently, $R$ is a left coherent ring. Hence, $\mathcal{FP}_n$-$Proj(R) = \mathcal{FP}_1$-$Proj(R)$ for all $n > 1$. So, by \cite[Proposition 3.7]{MD1}, $\mathcal{FP}_n$-$Proj(R)$ is closed under submodules.

$(3) \Rightarrow (4)$ 
Follows from \cite[Lemma 2.2]{PT}.

$(4) \Rightarrow (1)$ 
Follows from \cite[Theorem 3.16]{Shojaee}.
\end{proof}

As a consequence, if \(R\) is a left \(FP\)-hereditary ring, the class \(FP_n\text{-}Proj(R)\) is closed under direct products if and only if it is enveloping \cite[Lemma 2.2]{YT}.\\

Since $\npd_R(M)\leq\pdR(M)$ for any right $R$-module $M$, it is clear that $\gpd(R)\leq \glD(R)$, where $\glD(R)$ denotes the left global dimension of the ring $R$. For $n$-von Neumann regular rings, the equality holds. In fact, we have:

\begin{coro}\label{dim global Neumann}
Let $R$ be a left $n$-von Neumann regular ring. Then,
\begin{enumerate}
\item $\gpd(R)= \glD(R)$.
\item  $R$ is left hereditary if and only if every left ideal of $R$ is $FP_n$-projective.
\item  $R$ is left hereditary if and only if $R$ is left $FP$-hereditary.
\end{enumerate}
\end{coro}
\qed

For the general case, given a positive integer $k$, following  \cite[Theorem 4.4]{Stafforf-Warfiel} we can find a ring $R$  satisfying   $$\gpd(R) + k = \glD(R).$$ 

Now, we want to find the rings $R$ with $\gpd(R)= 1$. 

\begin{prop}
$\gpd(R)\leq 1$ if and only if $\npd_R(P/M)\leq 1$ for every quotient $P/M$ with $M$ a submodule of a projective $P$.
\end{prop}

\begin{proof}
We only need to prove the converse. Given any $R$-module $N$ we have a short exact sequence
$0\rightarrow \rm Ker(f) \rightarrow P \stackrel{f}{\rightarrow} N \rightarrow 0$ with $P$ projective. By assumption $\npd_R(N)=\npd_R(P/\rm Ker(f))\leq 1$ and equivalently $\gpd(R)\leq 1$.
\end{proof}

Given a class $\mathcal{A}$ of $R$-modules we denote by $\gpd(\mathcal{A})= \sup \{\fppd_R(M)\mid M \in \mathcal{A}\}$.\\

Observe that  $R$ is a strong left $(n-1)$-coherent ring if and only if $\gpd(\mathcal{FP}_{n-1}(R))=0$. In addition,  $R$ is a left $(n-1)$-coherent if and only if $\gpd(\mathcal{CFP}_{n-1}(R))=0$; where $\mathcal{CFP}_{n-1}(R)$ denotes  the class of all cyclic and  finitely $(n-1)$-presented $R$-modules.

Next proposition implies that in order to prove that a ring is strong left $(n-1)$-coherent it suffices to show that $\gpd(\mathcal{FP}_{n-1}(R)) \leq 1$.

\begin{prop}
Let $R$ be a ring and $n \geq 2$. If $\gpd(\mathcal{FP}_{n-1}(R)) \leq 1$, then $\gpd(\mathcal{FP}_{n-1}(R))=0$.
\end{prop}
\begin{proof}
Assume that there exists $M$ finitely $(n-1)$-presented with  $\fppd_R(M)= 1$. Then, for all $FP_n$-injective modules $E$, we have $\Ext^{2}_R(M, E) = 0$. Consider the short exact sequence $$0\rightarrow K \rightarrow P \rightarrow M \rightarrow 0,$$ where $P$ is projective and finitely generated and $K$ is finitely $(n-2)$-presented. This sequence induces an  exact sequence
$$\cdots\rightarrow 0=\Ext^{1}_R(P,E)\rightarrow \Ext^{1}_R(K,E) \rightarrow \Ext^{2}_R(M,E)=0 \rightarrow\cdots$$
which implies that $\Ext^{1}_R(K,E)= 0$. Therefore, $K \in \mathcal{FP}_{n}$-$Proj(R)\subseteq \mathcal{FP}_{n-1}$-$Proj(R)$.  Since $K$ is also finitely generated, by Lemma \ref{FP-proj lemma}, $K$ is in fact finitely $(n-1)$-presented. Consequently, $M$ is finitely $n$-presented. Finally, $\fppd_R(M) = 0$, contradicting our assumption.
\end{proof}

A direct consequence is the following corollary.

\begin{coro}\label{n-1 coherente}
A ring $R$ is strong left $(n-1)$-coherent if and only if $\gpd(\mathcal{FP}_{n-1}(R)) \leq 1$.
\end{coro}
\qed

Observe that the same conclusion can be drawn for $(n-1)$-coherent rings if we apply the argument of the previous proof again, with the class $\mathcal{FP}_{n-1}(R)$ replaced by the class $\mathcal{CFP}_{n-1}(R)$.

\begin{coro}
A ring $R$ is left $(n-1)$-coherent if and only if $\gpd(\mathcal{CFP}_{n-1}(R)) \leq 1$.
\end{coro}
\qed

If $n \geq 2$, $\gpd(\mathcal{FP}_{n-1}(R)) \leq 1$ implies $\gpd(\mathcal{CFP}_{n-1}(R)) \leq 1$. The converse holds for $n = 2$, i.e., every left coherent ring is strong left coherent. However, for $n > 2$, the converse is an open problem known as the $n$-coherence conjecture, i.e.,  it is not known if $(n-1)$-coherence implies strong $(n-1)$-coherence.\\

Our purpose now is to characterize the rings with $FP_n$-projective global dimension equal to 1.\\

\begin{obs}\label{FPn=fp si r es coherente}
If $R$ is a left coherent ring we have that $\mathcal{FP}_n$-$Inj(R)= \mathcal{FP}_1$-$Inj(R)$ for $n\geq 2$. Then, for every $R$-module $M$ we have $ \fppd_R(M)=\fpd(M)$. Therefore, $\gpd(R)=\fpD(R)$.
\end{obs}

\begin{prop}
Given $n\geq 2$, a ring $R$ is left $FP$-hereditary if and only if $\gpd(R)\leq 1$.
\end{prop}
\begin{proof} 
Assume that $R$ is left $FP$-hereditary. Then, $R$ is left coherent. By Remark \ref{FPn=fp si r es coherente}, for $n \geq 2$, we have $\gpd(R) = \mathrm{fpD}(R)$. Consequently, $\gpd(R) \leq 1$ by \cite[Proposition 3.7]{MD1}.

For the converse, assume that $R$ satisfies $\gpd(R) \leq 1$. Then, $\igpd(R)\leq 1$ for all $1 \leq i < n$. For $i=2$, Corollary \ref{n-1 coherente} gives that $R$ is a left coherent ring. Therefore, $ \fpD(R) =\unogpd(R) \leq 1$. Again,
\cite[Proposition 3.7]{MD1}, implies that $R$ is left $FP$-hereditary.
\end{proof}

\begin{obss} We note the following:
\begin{enumerate}
\item There exist rings for which the equality in the previous proposition holds. For example, let $\K$ be a field and $R = \K \langle x, y \rangle$ be the non-commutative polynomial ring in two variables. It can be shown that $R$ is an $FP$-hereditary ring but not a Noetherian ring. See \cite[Example 2.4.2]{Shojaee}.

\item It is always possible to construct a ring that is not strong left $(n-1)$-coherent \cite[Example 2]{Xue}. Therefore, Corollary \ref{n-1 coherente} implies that there exist rings for which $\gpd(R) \geq 2$.
\end{enumerate}
\end{obss}

 \subsection{$FP_n$-projective  dimension and  $\lambda$-dimension:}
 The \textit{$\lambda$-dimension} of an $R$-module $M$, denoted by $\lamd_R(M)$, was defined  in \cite{Bou} as follows:
\[ \lamd_R(M) \left\{
\begin{array}{lll}
= \infty & \text{if } & M\in \mathcal{FP}_{\infty}(R).\\
= n & \text{if } & M\in \mathcal{FP}_n(R)\setminus \mathcal{FP}_{n+1}(R).\\
\geq n &\text{if } & M\in \mathcal{FP}_n(R).\\
= -1 &\text{if } & M \in \RMod \setminus \mathcal{FP}_0(R).
\end{array}
\right.
\]
In addition, the $\lambda$-dimension of the ring $R$, denoted by $\lambda$-$\dim(R)$, is defined as the least integer $n$ (or $\infty$ if no such integer exists) such that $\lambda$-$\dim_R(M) \geq n$ implies $\lambda$-$\dim_R(M) = \infty$. It is well-known that $R$ is left Noetherian if and only if $\lambda$-$\dim(R)=0$, and $R$ is strong left $n$-coherent if and only if $\lambda$-$\dim(R) \leq n$.

The $\lambda$-dimension and the global $FP_n$-projective dimension  are distinct and can have different values for a given ring $R$. We can make the following considerations:

\begin{enumerate}
\item $R$ is left Noetherian if and only if $\lambda$-$\dim(R) = \gpd(R)= 0$.
\item If $\gpd(R) \leq 1$, it follows that $\lambda$-$\dim(R) \leq 1$ because every left $FP$-hereditary ring is left coherent. However, \cite[Example 3.9]{Shojaee} shows a ring $R$ for which $\lambda$-$\dim(R) \leq 1$ and $\gpd(R) \geq 2$.   Furthermore, \cite[Example 2.4]{Mahdou3} presents a local 2-von Neumann regular ring $R$ with $\glD(R)=\infty$.  Therefore,  $\lambda$-$\dim(R) \leqslant 2$ by \cite[Theorem 4.8]{MD4}  and by Corollary \ref{dim global Neumann}, $\gpd(R) = \infty$.
\end{enumerate}

\subsection{The $FP_n$-projective dimension over strong $n$-coherent rings.}

As previously mentioned, the $FP_1$-projective dimension of a module coincides with the $FP$-projective dimension $\rm fpd$ defined and studied in \cite{MD1}. Moreover, over left coherent rings, the $FP_1$-projective global dimension agrees with the $FP$-projective global dimension $\rm fpD$ (see \cite[Theorem 3.1]{MD1}). Therefore, below we will focus on the case $n \geq 2$.

The $FP_n$-projective dimension of an $R$-module coincides with the $(n,0)$-projective dimension defined and studied in \cite{Ouy}.  Over strong left $n$-coherent rings the $FP_n$-projective dimension of an $R$-module $M$ also coincides with the relative projective dimension $\rm pd_{\mathcal{X}}(M)$ of $M$ with respect to the class  $\mathcal{X}=\mathcal{FP}_n$-$Inj(R)$, defined by Auslander and Buchweitz in \cite{AB}.  Over strong left $n$-coherent rings, the $FP_n$-projective global dimension of $R$ agrees with the $(n,0)$-projective dimension of $R$ (see \cite[Theorem 3.1]{Ouy}) and with the relative projective dimension $\rm pd_{\mathcal{X}}(R\text{-}Mod)$; see \cite{AM}.

Moreover, following \cite[Proposition 3.1]{Ouy}, we can see that over strong left $n$-coherent rings, the $FP_n$-projective global dimension of $R$ agrees with the left global projective dimension relative to
$\mathcal{X}= \mathcal{FP}_n$-$Proj(R)$ denoted by $\text{PD}_{\mathcal{X}}(R)$ in the sense of \cite[Definition 3.1]{Cortez}. Observe that, in \cite{Cortez},  $\mathcal{X}$ is a class of $R$-modules containing all projective modules and the projective dimension relative to $\mathcal{X}$ of a module $M$ is defined using resolutions (i.e., $M$ has projective dimension relative to $\mathcal{X}$, or $\mathcal{X}$-projective dimension less than or equal to $k$, if there exists a projective resolution of $M$ such that its $(k-1)^{\text{st}}$ syzygy belongs to $\mathcal{X}$). 

As a consequence of \cite[Lemma 3.9]{Cortez} taking $\mathcal{X}= \mathcal{FP}_n$-$Proj(R)$ and $\mathcal{Y}= \text{R-Mod}$ we obtain the following corollary.

\begin{prop}
Let $R$ be a strong left $n$-coherent ring. Then, each $R$-module has a finite $FP_n$-projective dimension if and only if $\gpd(R) < \infty$.
\end{prop}
\qed

We state the following theorem which can  be deduced from \cite[Lemma 1.1]{AM} and \cite[Corollary 2.3]{AM} for convenient reference. 

\begin{teo}  \label{Fp-injec} Let $R$ be a strong left  $n$-coherent ring. Then the following are identical.
\begin{enumerate}
\item[(1)] $\gpd(R)$.
\item[(2)] $\sup \{\fppd_R(M)\mid M \text{ is an $FP_n$-injective $R$-module} \}$.
\item[(3)] $\sup \{\id_R(M) \mid M \text{ is an $FP_n$-injective $R$-module} \}$.
\item[(4)] $\sup \{\fppd_R(M)\mid M \text{ is a cyclic $R$-module} \}$.
\item[(5)] $\sup \{\fppd_R(M)\mid M \text{ is a finitely generated $R$-module} \}$.
\end{enumerate}
\end{teo}
\qed

The following result is a fairly straightforward of Theorem  \ref{Fp-injec} and the fact that a $R$-module $M$ is $FP_n$-projective if and only if $\fppd_R(M)=0$.

\begin{coro}\label{coro fp proyective dimension} Let $R$ be a strong left $n$-coherent ring. An  $R$-module $M$ is $FP_n$-projective if and only if  $\Ext^{j}_R(M,N)=0$ for all $FP_n$-injective  $R$-module $N$ and all positive integer $j$.
\end{coro}
\qed

\begin{defi} Let \( R \) be a ring. The \textit{finitistic \( FP_n \)-projective global dimension of \( R \)} is defined as
 $$\fgpd(R)=\sup \{\fppd_R(M)\mid M \text{ is an $R$-module with} \fppd_R(M) < \infty \}.$$
\end{defi}

Clearly, if every $R$-module has a finite $FP_n$-projective dimension, then $\fgpd(R) = \gpd(R)$.

\begin{coro}
Let $R$ be a strong left $n$-coherent ring such that every $R$-module has finite $FP_n$-projective dimension. Then,  $\fgpd(R) < \infty$.
\end{coro}
\qed

\begin{prop}
Let $R$ be a strong left $n$-coherent ring. The following statements are equivalent. 
\begin{enumerate}
    \item $\fgpd(R)=0$.
    \item If $P$ and $F$ are $FP_n$-projective $R$-modules with $P \subseteq F$, then $F/P$ is $FP_n$-projective.
    \item If $M$ has a finite $FP_n$-projective resolution, then $M$ is $FP_n$-projective.
\end{enumerate}
\end{prop}

\begin{proof}
    $(1) \Leftrightarrow (3)$ Follows from Theorem \ref{Fp-injec}.

    $(1) \Rightarrow (2)$ Suppose that $\fgpd(R) = 0$ and let $P$ and $F$ be $FP_n$-projective $R$-modules with $P \subseteq F$. Since the sequence $0 \rightarrow P \rightarrow F \rightarrow F/P \rightarrow 0$ is exact, $F/P$ has a finite $FP_n$-projective resolution and then  $\fppd_R(F/P) < \infty$. Therefore, $\fppd_R(F/P) = 0$, meaning that $F/P$ is $FP_n$-projective.

$(2) \Rightarrow (1)$ Let $M$ be an $R$-module with $\fppd_R(M) = k < \infty$. Then $M$ has a finite $FP_n$-projective resolution:
\[
0 \rightarrow P_k \rightarrow P_{k-1} \rightarrow \cdots \rightarrow P_1 \rightarrow P_0 \rightarrow M \rightarrow 0,
\]
where each $P_i$ is $FP_n$-projective for every $0 \leq i \leq k$. Let $N_{k-2}$ denote the $(k-2)$-th syzygy. If $k > 0$, then by hypothesis, $N_{k-2} \cong P_{k-1}/P_k$ is $FP_n$-projective. From the exact sequence
\[
0 \rightarrow N_{k-2} \rightarrow P_{k-2} \rightarrow \cdots \rightarrow P_1 \rightarrow P_0 \rightarrow M \rightarrow 0,
\]
it follows that $\fppd_R(M) \leq k - 1$, which contradicts the assumption that $\fppd_R(M) = k$. Thus, $k = 0$, meaning that $M$ is $FP_n$-projective.
\end{proof}

\section{Weak and projective dimension of $FP_n$-projective modules}

In \cite{Zhu} Zhu introduced the following  dimension of an $R$-module $M$.

\begin{defi} For any left  (or right) $R$-module $N$, the \textit{$FP_n$-flat dimension} of $N$, denoted by $\nfd_{R}(N)$, is the smallest integer $k \geq 0$ such that $\Tor^{R}_{k+1}(F, N) = 0$ (or $\Tor^{R}_{k+1}(N, F) = 0$, respectively) for every finitely $n$-presented right (or left, respectively) $R$-module $F$. If no such integer $k$ exists, we define $\nfd_{R}(N) = \infty$.
\end{defi}

\begin{obs}\label{fpn-flat}
By \cite[Theorem 2.3.6]{ZP}, $\nfd_{R}(N)$ is the smallest non-negative integer $k$ such that $N$ has a resolution by $FP_n$-flat modules:
\[
\cdots \rightarrow Q_i \rightarrow Q_{i-1} \rightarrow \cdots \rightarrow Q_1 \rightarrow Q_0 \rightarrow N \rightarrow 0
\]
where $Q_i = 0$ for every $i > k$. Furthermore, $\nfd_{R}(N) \leq k$ if and only if every projective $k$-th syzygy of $N$ is $FP_n$-flat.
\end{obs}

For a left coherent ring $R$, it is known that the projective dimension $\pdR(M)$ of any finitely $1$-presented $R$-module $M$ is equal to its weak dimension $\wdR(M)$. Moreover, if in addition $R$ is left self $FP_1$-injective, the equality of both dimensions holds true for any $FP_1$-projective $R$-module  $M$  \cite[Proposition 4.1]{MD1}.

In the following proposition we have a generalization  for $n\geq 2$.

\begin{prop}\label{coro wd=pd} Let $R$ be a strong left $n$-coherent and left self $FP_n$-injective ring. If $M$ is an $FP_n$-projective  $R$-module, then $\pdR(M)= \nfd_{R}(M)$.
\end{prop}

\begin{proof} 
Clearly, for every $FP_n$-projective $R$-module $M$, we have $\nfd_{R}(M) \leq \wdR(M) \leq \pdR(M)$.

Assume that $\nfd_{R}(M) = k < \infty$. By Remark \ref{fpn-flat}, there is an exact sequence of $R$-modules:
\[
0 \rightarrow F_k \rightarrow P_{k-1} \rightarrow \cdots \rightarrow P_1 \rightarrow P_0 \rightarrow M \rightarrow 0
\]
where $P_i$ is projective for $0\leq i \leq k-1$ and $F_k$ is $FP_n$-flat. It is enough to show that $F_k$ is also projective.

Now, there is a short exact sequence of $R$-modules:
\[
0 \rightarrow K \rightarrow F \rightarrow F_k \rightarrow 0
\]
where $F$ is a free module and by hypothesis, $F$ is $FP_n$-injective. By Lemma \ref{pur2}, this short exact sequence is $n$-pure.

Let $N$ be an $FP_n$-injective $R$-module. By Corollary \ref{coro fp proyective dimension}, we have $\Ext^{1}_R(F_k, N) \simeq \Ext^{k+1}_R(M, N) = 0$. This implies that $F_k$ is $FP_n$-projective and hence a direct summand of $F$ by Theorem \ref{FP-proj pureza}. Therefore, $F_k$ is projective.
\end{proof}

\begin{coro}  Let $R$ be a strong left $n$-coherent and left self $FP_n$-injective ring. If $M$ is a non-projective  $FP_n$-projective $R$-module, then $\nfd_{R}(M)=\infty$.
\end{coro}

\begin{proof}
 By Proposition \ref{coro wd=pd}, it is enough to show that $\pdR(M)=\infty$. Suppose that $\pdR(M)= k < \infty$. Since $M$ is non-projective $k\geq 1$. Then, there exists an $R$-module $N$ such that $\Ext_R^k(M,N) \neq 0$. Consider the short exact sequence of $R$-modules $0 \rightarrow L \rightarrow P \rightarrow N \rightarrow 0$ where $P$ is projective. We consider the induced exact sequence

$$\Ext^{k}_R(M,P)\rightarrow \Ext^{k}_R(M,N) \rightarrow \Ext^{k+1}_R(M,L)=0. $$

Since $R$ is left self $FP_n$-injective, $P$ is also $FP_n$-injective. Then, $\Ext_R^k(M,P)= 0$  by Corollary \ref{coro fp proyective dimension}. It follows that $\Ext_R^k(M,N) =0$, a contradiction.
\end{proof}

Observing that every finitely $n$-presented $R$-module is $FP_n$-projective, we have the following.

\begin{coro}  Let $R$ be a strong left $n$-coherent and left self $FP_n$-injective ring. If $M$ is a non-projective  finitely $n$-presented $R$-module, then $\nfd_{R}(M)=\infty$.
\end{coro}
\qed

The following corollary is an immediate consequence of the previous results and the fact that for every $R$-module $M$ we have $\nfd_{R}(M) \leq \wdR(M) \leq \pdR(M)$.

\begin{coro}\label{coro debil}  Let $R$ be a strong left $n$-coherent and left self $FP_n$-injective ring. 
\begin{enumerate}
\item If $M$ is an $FP_n$-projective  $R$-module, then $\pdR(M)= \wdR(M)$.
\item If $M$ is a non-projective $FP_n$-projective  $R$-module, then $\wdR(M)=\infty$.
\item If $M$ is a non-projective  finitely $n$-presented $R$-module, then $\wdR(M)=\infty$.
\end{enumerate}
\end{coro}
\qed

Given an integer $k\geq 1$, recall that an $R$-module $T$ is  said to be \textit{$k$-tilting} when it satisfies:
\begin{itemize}
    \item[(T1)] $\pdR(T)\leq k$
    \item[(T2)] $\Ext^i_R(T, T^{(I)}) = 0$ for each $i \geq 1$ and all sets $I$, and
    \item[(T3)] there exist $r \geq 0$ and a long exact sequence
    $0 \rightarrow R \rightarrow T_0 \rightarrow \cdots \rightarrow T_r \rightarrow 0 $ such that $T_i \in \text{Add}(T)$ for each $0 \leq i \leq r$.\\
\end{itemize}

Here, $\text{Add}(T)$ denotes the class of all direct summands of arbitrary direct sums of copies of $T$. \\

A class of $R$-modules $\mathcal{X}$ is $k$-tilting if there is a $k$-tilting $R$-module $X$ such that $\mathcal{X} = X^{\perp}$. A cotorsion pair $(\mathcal{A}, \mathcal{B})$ is $k$-tilting provided that $\mathcal{B}$ is a $k$-tilting class.

An immediate consequence of the previous corollary and \cite[Lemma 1.13]{PT} is the following:

\begin{coro}
    Let $R$ be a strong left $n$-coherent and left self $FP_n$-injective ring. The hereditary cotorsion pair $(\mathcal{FP}_n$-$Proj(R), \mathcal{FP}_n$-$Inj(R))$ is not $k$-tilting for any $k \geq 1$.
\end{coro}
\qed

Clearly every left \( n \)-von Neumann regular ring is strong left \( n \)-coherent and left self-\( FP_n \)-injective. Moreover, 
if we assume that $R$ is a strong $n$-coherent domain, then $R$ is self $FP_n$-injective if and only if $R$ is a field  \cite[Proposition 2.1]{Xing}. The following example shows the existence of a ring which is strong $2$-coherent  and self $FP_2$-injective, but not  coherent.

\begin{ejem} Let $R$ be a local Noetherian regular ring with Krull dimension equal $n$, with $\mathsf{m}$ being its maximal ideal and $E(R/\mathsf{m})$ being the $R$-injective envelope of  $R/\mathsf{m}$. If $R$ is complete in its $\mathsf{m}$-adic topology, then the  trivial extension (denoted by $R\ltimes E(R/\mathsf{m})$) of $R$ by $E(R/\mathsf{m})$ is a strong $n$-coherent ring \cite[Theorem A']{Roos} and self $FP_n$-injective ring \cite[Remark II.6]{Couchot3}. However, it is not a $(n-1)$-coherent ring.
\end{ejem}

Another known example of  a ring which is strong $2$-coherent  and self $FP_2$-injective but not  coherent, is the polynomial ring 
${\displaystyle R = \dfrac{\mathsf{k}[x_1, x_2, \dots]}{(x_i x_j)_{i,j \geq 1}}}$  with \( \mathsf{k} \) a field. See \cite[Examples 1.3 and 5.7]{BP}.

\subsection{Applications to valuation rings}

Let $R$ be a commutative ring. Recall that an $R$-module $M$ is said to be \textit{uniserial} if the set of its submodules is totally ordered by inclusion. A ring $R$ is a \textit{valuation ring}  if it is uniserial as $R$-module. We note that $R$ is a valuation ring if and only if $R$ is a local ring and all finitely generated ideals are principal.\\

Valuation rings are always strong $2$-coherent and are coherent if the ring does not have nonzero zero divisors \cite[Theorem II.11]{Couchot3}. Therefore,  for every self $FP_2$-injective valuation ring,   Corollary \ref{coro debil} implies that $\pdR(M)= \wdR(M)$ for all $FP_2$-projective $R$-module $M$.

 For example, if $R$ is a valuation ring with maximal ideal $\mathfrak{m}$ equal to the set of zero divisors,  then $R$ is strong $2$-coherent and self $FP_1$-injective (hence self $FP_2$-injective) by \cite[Theorem II.11]{Couchot3}.

A ring $R$ is \textit{arithmetical} if it is locally a valuation ring.
Every arithmetical ring is strong $3$-coherent \cite[Theorem II.1]{Couchot3}. According to \cite[Theorem 1]{Couchot12}, any module over an arithmetical ring $R$ has weak dimension $0, 1, 2$ or $\infty$. If $R$ is also self $FP_3$-injective, Corollary \ref{coro debil} implies that $\pdR(M)= 0, 1, 2$ or $\infty$ for every $FP_3$-projective $R$-module $M$.

Moreover, we can get similar results for other arithmetical rings, ensuring that they are  coherent or strong $2$-coherent.  For example, if $R$ is also self $FP_1$-injective, then $R$ is strong $2$-coherent \cite[Theorem II.4]{Couchot3}. Additionally, if the annihilator of every element in $R$ is finitely generated, $R$ is coherent \cite[1.4 Fact, Ch XII, \textsection Arithmetic Rings ]{Lombardi}. On the other hand, if $R$ is also reduced, then $R$ is strong $2$-coherent and it is well known that $R$ has global weak dimension at most $1$. We summarize this in the following corollary.

\begin{coro} Let $R$ be an arithmetical ring.
\begin{enumerate}
\item If $R$ is self $FP_1$-injective ring, then $\pdR(M)= 0, 1, 2$ or $\infty$ for every $FP_2$-projective $R$-module $M$.
\item If $R$ is self $FP_1$-injective ring and the annihilator of every element is finitely generated, then $\pdR(M)= 0, 1, 2$ or $\infty$ for every $FP_1$-projective $R$-module $M$.
\item  If $R$ es reduced and self $FP_2$-injective ring, then $\pdR(M)= 0$ or $1$ for every $FP_2$-projective $R$-module $M$.
\end{enumerate}
\end{coro}
\qed

Every arithmetical ring of Krull dimension zero is strong $2$-coherent and self $FP_2$-injective. See \cite[Corollary II.7 ]{Couchot3}.

\section{Relation with the  $FP_n$-injective  dimension of a ring}

In \cite{Zhu}, Zhu introduces the following dimension of an $R$-module $M$ and the respective global dimension of the ring $R$.

\begin{defi}
Let be $M$ an $R$-module. The \textit{$FP_n$-injective dimension} of $M$, which we denote by $\nid_R(M)$, is given by the smallest integer $k\geq 0$ such that $\Ext^{k+1}_R(F,M)=0$ for every $F\in \mathcal{FP}_n(R)$. If no such integer $k$ exists, we define $\nid_R(M) = \infty$.
\end{defi}

Recall that the \textit{left $FP_n$-injective global dimension} of a ring $R$, which we denote by $\gid(R)$, is defined by
$$\gid(R)= \sup \{\nid_R(M)\mid M \text{ is an $R$-module} \}.$$

\begin{obs}\label{dim fp-inj y sec} 
Over strong left $n$-coherent rings the $FP_n$-injective dimension of an 
$R$-module $M$ coincides with its relative injective dimension 
$\id_{\mathcal{X}}(M)$ with respect to 
$\mathcal{X}=\mathcal{FP}_n(R)$, and the $FP_n$-injective global dimension of 
$R$ agrees with $\id_{\mathcal{X}}(R\text{-}Mod)$; see \cite{AM}. Following \cite{AB}, if $\mathcal{X}$ is a class of $R$-modules, the $\mathcal{X}$-coresolution dimension $coresdim_{\mathcal{X}}(M)$  of $M$ is the minimal nonnegative integer $k$ such that there is an exact sequence $0 \to M \to X_0 \to \cdots \to X_k \to 0,\ 
X_i$ with $X_i\in \mathcal{X}$ for $0\le i \le k$.  With $coresdim_{\mathcal{X}}(M) := \infty$ if no such $k$ exists.  For a class $\mathcal{Y}\subseteq R\text{-}Mod$, we set
\[
coresdim_{\mathcal{X}}(\mathcal{Y})
:= \sup\{\, coresdim_{\mathcal{X}}(Y) : Y \in \mathcal{Y}\,\}.
\]

By \cite[Theorem~2.12]{Zhu}, for 
$\mathcal{X}=\mathcal{FP}_n\text{-}Inj(R)$ and 
$\mathcal{Y}=R\text{-}Mod$, we have
\[
coresdim_{\mathcal{X}}(\mathcal{Y}) = \gid(R).
\]
\end{obs}

\begin{prop}\label{prop relacion dimensiones} Let $R$ be a strong left $n$-coherent ring. Then the following are identical.
\begin{enumerate}
\item  $ \gid(R)$.
\item  $\sup \{\pdR(M)\mid M \text{ is an $FP_n$-projective  $R$-module} \}.$
\item  $\sup \{\pdR(M)\mid M \text{ is a finitely $n$-presented $R$-module} \}.$
\item $\sup \{\nid_R(M)\mid M \text{ is an $FP_n$-projective  $R$-module} \}.$
\item $\sup \{\nfd_{R}(M)\mid M \text{ is a right $R$-module} \}.$
\end{enumerate}
\end{prop}
\begin{proof}

$(1)=(5)$ Follow from \cite[Theorem 3.8]{Gao}.\\
$(1)=(3)$ Follow from \cite[Lemma 1.1]{AM}. \\
$(1)=(2)=(4)$ Follow from \cite[Theorem 2.5]{AM} and Remark \ref{dim fp-inj y sec}.

\end{proof}

Since $R$ is  left $n$-von Neumann regular ring if and only if  every $R$-module is $FP_n$-injective by \cite[Theorem 3.9]{Zhu} or equivalently $\gid(R)= 0$, we recover part of \cite[Theorem 4.8]{MD4}.

\begin{coro} Let $R$ be a ring. Then the following conditions are equivalent.
\begin{enumerate}
\item $R$ is  left $n$-von Neumann regular ring.
\item $\gid(R)= 0$.
\item $R$ is strong left $n$-coherent ring and every $FP_n$-projective $R$-module is $FP_n$-injective.
\end{enumerate}
\end{coro}
\qed

We recall that given two integers $n,d\geq 0$, a ring $R$ is said to be a left \textit{$(n,d)$-ring} if every finitely $n$-presented  module has projective dimension at most $d$ \cite{Costa}. Observe that when $d = 0$, one gets back
the left $n$-von Neumann regular rings and if $d = 1$, one obtains
the left $n$-hereditary rings.

\begin{coro}
Let $R$ be a strong left $n$-coherent ring and $k \geq 0$ an integer. Then $R$ is a left $(n,k)$-ring if and only if $\gid(R) \leq k$.
\end{coro}
\qed

Following \cite[Theorem 4.5]{Costa}, if \( R \) is a ring with weak dimension less than or equal to \( k \), then \( R \) is a \((k+1, k)\)-ring. Moreover, it is known that these rings are strong \( k+1 \)-coherent. \cite[Theorem 2.2]{Costa}

\begin{coro}
   Let \( R \) be a ring and let \( n \) and \( k \) be non-negative integers with \( n \leq k \). If \( R \) is strong left \( n \)-coherent with weak finite dimension equal to \( k \), then \( R \) is a left \((n, k)\)-ring.
\end{coro}
\qed

In \cite{BPa} $n$-hereditary rings are characterized by  bounding  by $1$ the projective (or weak) dimension over the class of its finitely $n$-presented modules. We show that the same holds true if we change the class of finitely $n$-presented modules for the class of  $FP_n$-projective modules.

\begin{prop}\label{n-hereditario} Let $R$ be a ring. The following conditions are equivalent.
\begin{enumerate}
\item $R$ is a left $n$-hereditary ring.
\item $\pdR(M)\leq 1$  for every $FP_n$-projective $R$-module $M$.
\item $\wdR(M)\leq 1$  for every $FP_n$-projective $R$-module $M$.
\end{enumerate}
\end{prop}
\begin{proof}
$(1)\Rightarrow(2)$ Follows from \cite[Theorem 3.2]{Zhu} and Proposition \ref{prop relacion dimensiones}.

$(2)\Rightarrow(3)$ It is clear.

$(3)\Rightarrow(1)$ Follows from  \cite[Lemma 3.2]{BPa}.
\end{proof}

If \( n \geq 2 \) and \(\mathrm{w.dim}(R) \leq 1\), then \( R \) is both left and right \( n \)-hereditary according to \cite[Corollary 3.6]{Zhu}. By applying Proposition \ref{prop relacion dimensiones} and considering the \( FP_n \)-flat dimension, we can recover \cite[Corollary 3.5]{Zhu2}.

\begin{prop}\cite[Corollary 3.5]{Zhu2}  The following conditions are equivalent for a ring $R$.
\begin{enumerate}
\item $R$ is a left $n$-hereditary ring.
\item Every submodule of an $FP_n$-flat right $R$-module is $FP_n$-flat.
\item Every right ideal of $R$ is $FP_n$-flat.
\item Every finitely generated right ideal of $R$ is $FP_n$-flat.
\end{enumerate}
\end{prop}
\qed

\section{Further consequences}

We finish, for $n \geq 2$, with some applications to subprojectivity domains, the CF-conjecture and trace modules

\subsection{Subprojectivity domains.}
We now present the relationship between \( FP_n \)-projective modules and subprojectivity domains over strong left \( n \)-coherent and left self \( FP_n \)-injective rings. Following \cite{Holston}, given two \( R \)-modules \( M \) and \( N \), \( M \) is said to be \textit{\( N \)-subprojective} if for every epimorphism \( g : B \to N \) and for every homomorphism \( f : M \to N \), there exists a homomorphism \( h : M \to B \) such that \( gh = f \). For an \( R \)-module \( M \), the \textit{subprojectivity domain} of \( M \) is defined to be the collection of all \( R \)-modules \( N \) such that \( M \) is \( N \)-subprojective. We denote the subprojectivity domain of \( M \) by \(\text{Proj}^{-1}(M)\).

The subprojectivity domain, or domain of subprojectivity, of a class of \( R \)-modules \(\mathcal{M}\) is defined as
\[ \text{Proj}^{-1}(\mathcal{M}) := \{ N\in \rm R\text{-}Mod \mid M \text{ is } N\text{-subprojective for every } M \in \mathcal{M} \}. \]

Clearly every module in \(\mathcal{M}\) is projective if and only if \(\text{Proj}^{-1}(\mathcal{M}) = \mathrm{R\text{-}Mod}\). \\

Assuming that $\mathcal{M} = \mathcal{FP}_n$-$Proj(R)$ and that all projective $R$-modules belong to $\mathcal{FP}_n$-$Inj(R)$ (for example when \( R \)  is a left self \( FP_n \)-injective ring) then, following the comments before \cite[Theorem 4.1]{Amzil}, we obtain:

\begin{prop}\label{prop7}
    Let \( R \) be a left self \( FP_n \)-injective ring and \( N \) be an \( R \)-module. Then, 
    $$ N \in \text{Proj}^{-1}(\mathcal{FP}_n\text{-Proj}(R)) \text{ if and only if } \Omega^1_R(N) \in \mathcal{FP}_n\text{-Inj}(R)  $$
    
    where $\Omega^1_R(N)$ denotes the first syzygy of \( N \).
\end{prop}
\qed

Recall that an \( R \)-module \( M \) is called \textit{subprojectively poor} (or \textit{sp-poor} or \textit{p-indigent}) if its subprojectivity domain consists of only projective modules. A natural question to ask is how small \(\text{Proj}^{-1}(M)\) can be. \cite[Proposition 2.8]{Holston} shows that the domain of subprojectivity of any module must contain at least the projective modules.

\begin{prop}
    Let \( R \) be a left self \( FP_n \)-injective ring. If \(\text{Proj}^{-1}(\mathcal{FP}_n\text{-Proj}(R)) \subseteq \mathcal{FP}_n\text{-Proj}(R)\), then \(\text{Proj}^{-1}(\mathcal{FP}_n\text{-Proj}(R)) = \text{Proj}(R)\)
\end{prop}

\begin{proof}
 It is sufficient to show that $\text{Proj}^{-1}(\mathcal{FP}_n$-$Proj(R)) \subseteq \text{Proj}(R)$. Let \( N \) be an \( R \)-module such that $ N \in \text{Proj}^{-1}(\mathcal{FP}_n$-$Proj(R))$. By Proposition \ref{prop7}, $\Omega^1_R(N) \in \mathcal{FP}_n$-$Inj(R)$ and  by hypothesis $N\in \mathcal{FP}_n$-$Proj(R)$, then  \( \Ext^{1}_R(N, \Omega^1_R(N)) = 0 \). Therefore, \( N \) is projective by \cite[Lemma 1.1]{Trlifaj3}.
\end{proof}

Recall that, given a class $\mathcal{F}$ of $R$-modules, an \textit{$\mathcal{F}$-precover} of an $R$-module $M$ is a homomorphism $F \to M$ with $F \in \mathcal{F}$, such that 
\[
\text{Hom}_R(F', F) \to \text{Hom}_R(F', M) \to 0 
\]
is exact for any $F' \in \mathcal{F}$. An $\mathcal{F}$-precover is said to be \textit{special} provided that it is an epimorphism with kernel in the class $\mathcal{F}^\perp$. $\mathcal{F}$-preenvelopes and special $\mathcal{F}$-preenvelopes are defined dually.

\begin{prop}
    Let \( R \) be a strong left \( n \)-coherent ring. Then the following statements are equivalent.
    \begin{enumerate}
        \item \( \mathcal{FP}_n\text{-Inj}(R)= \text{Proj}^{-1}(\mathcal{FP}_n\text{-Proj}(R))\).
        \item \(\mathcal{FP}_n\text{-Proj}(R) \cap \mathcal{FP}_n\text{-Inj}(R) = \text{Proj}(R)\) and every \( FP_n \)-injective \( R \)-module  has a special \(\mathcal{FP}_n\text{-Proj}(R)\)-precover.
        \item \( R \) is a left self \( FP_n \)-injective, \(\text{Proj}^{-1}(\mathcal{FP}_n\text{-Proj}(R))\) is closed under cokernels of monomorphisms and every \( FP_n \)-projective \( R \)-module has a \(\mathcal{FP}_n\text{-Inj}(R)\)-preenvelope which is projective.
    \end{enumerate}
\end{prop}
\begin{proof}
     By \cite[Theorem 4.1]{MD4}, a ring \( R \) is strong left \( n \)-coherent if and only if \( \mathcal{FP}_n$-$Proj(R)$ is closed under kernels of epimorphisms. Therefore, the result follows from \cite[Theorem 4.1]{Amzil}.
\end{proof}

\subsection{Applications to $CF$-conjecture}
Rings satisfying that every cyclic module can be  embedded  in a free module are usually called in the literature left \textit{$CF$-rings}. The question of whether  any left $CF$-ring is left Artinian is nowadays known as the $CF$-conjecture. This conjecture is still open in general but it is known to be true under many different additional hypothesis, for example, being in addition  a perfect ring. See \cite[Theorem 3.4]{RS}.

According to \cite{Azumaya}, an $R$-module $M$ is called  \textit{finitely projective} if, for any finitely generated submodule $N$, the inclusion map $N \rightarrow M$ factors through a free $R$-module. It is well-known that finitely projective modules are always flat and if the ring $R$ is left Noetherian, then the converse is also true. See \cite{Azumaya}. 

By \cite[Proposition 3.3]{Couchot1}, if $R$ is a left self $FP_1$-injective ring in which every flat $R$-module is finitely projective, then $R$ is left perfect. This is the case of the rings  where every flat module is finitely projective and almost $FP_1$-injective, see  Remark \ref{almos}. This gives a new case where the $CF$-conjecture holds.

\begin{prop}\label{Azu-Couch}
Let $R$ be a left $CF$-ring satisfying that every flat $R$-module is almost $FP_1$-injective and finitely projective. Then $R$ is  left Artinian.
\end{prop}
\qed

It is clear that every left coherent ring is also a strong left $2$-coherent ring. We will show that if in addition $R$ is a left $CF$-ring, the converse holds true.  Therefore, we can extend the well known characterization of $CF$-rings to strong $2$-coherent rings.

\begin{prop}\label{CF} Let $R$ be a left $CF$-ring. The following are equivalent.
\begin{enumerate}
\item $R$ is a left Noetherian ring.
\item $R$ is a left coherent ring.
\item $R$ is a strong left $2$-coherent ring.
\item $R$ is a left $FP$-hereditary ring.
\end{enumerate}
\end{prop}

\begin{proof}
Since the equivalence of $(1)$, $(2)$, and $(4)$ was established in \cite[Corollary 3.10]{Shojaee}, and it is clear that $(2)$ implies $(3)$, it suffices to show that $(3)$ implies $(2)$. Assume that $R$ is a strong left $2$-coherent ring and let $I$ be a finitely generated ideal of $R$. We will show that $I$ is finitely $1$-presented. Since $R/I$ is cyclic, by hypothesis, it can be embedded in a free $R$-module $F$. Then, by \cite[Theorem 2.1]{Zhu}, $R/I$ is finitely $2$-presented and hence $I$ is finitely $1$-presented.
\end{proof}

\begin{obs}
If $R$ is a strong left $(n+1)$-coherent and $CF$-ring, then it is also left $n$-coherent. The proof follows a similar line of argument as the proof of ($3$) implies ($2$) in the previous proposition.
\end{obs}

Recall that $R$ is called a left \textit{Kasch ring} if every simple module embeds in $R$, or equivalently, if every simple $R$-module embeds in a free module. Artinian rings are always Kasch rings.

An immediate consequence of Proposition \ref{CF} and \cite[Lemma 2.10]{Garkusha} is that every strong left $2$-coherent and $CF$-ring $R$ is a left Kasch ring.

In order to prove that it is also an Artinian ring, we need an additional condition over the ring $R$. In fact, summarizing several known results, we have the following.

\begin{prop}\label{prop-cf-comm} Every strong left $2$-coherent $CF$-ring $R$ that satisfies at least one of the following conditions is Artinian.
\begin{enumerate}
\item Every flat $R$-module is almost $FP_1$-injective.
\item $R$ is a left or right semi-Artinian ring.
\item $R$ is a semiregular ring.
\item $R$ is a semiperfect ring.
\item The socle of $R$ is an essential submodule in the module $R$.
\item $R$ is a self-injective ring.
\item Every cyclic $R$-module is pure-injective (i.e., is injective with respect to pure exact sequences).
\item $R$ is left pure semisimple (i.e., if every $R$-module is pure-injective).

\end{enumerate}
\end{prop}
\begin{proof}
 $(1)$ By Proposition \ref{CF}, $R$ is a left Noetherian ring, so every flat $R$-module is finitely projective. Therefore, by assumption, every flat $R$-module is almost $FP_1$-injective and flat. Applying Proposition \ref{Azu-Couch}, it follows that $R$ is left Artinian.

 $(2),(3),(4)$ and $(5)$  Follow from Proposition \ref{CF} and \cite[Lemma 2.13]{Garkusha}.

 $(6)$ Follows from Proposition \ref{CF}, \cite[Prop 3.13]{Shojaee} and $(4)$.

 $(7)$ Follows from Proposition \ref{CF}, \cite[Lemma 2]{Ali} and $(4)$.

 $(8)$ Follows from $(7)$.

\end{proof}

\subsection{Trace modules in $FP_n$-injective envelopes}

 For basic terminology using here we refer  to \cite{Lindo}. Let $R$ be a ring and let $M$ and $X$ be $R$-modules. The \textit{trace module of $M$} in $X$ is the sum of all $R$-homomorphic images of $M$ in $X$, and it is denoted by $\tau_M(X)$:

$$\tau_M(X)=\sum_{\alpha\in Hom_R(M,X)} \alpha(M).$$

If there exists an injection $\iota: M \rightarrow X$ such that $\text{Img}(\iota)$ is a trace submodule of $X$, then $M$ is said to be \textit{trace in $X$ up to isomorphism}. As noted by Lindo and Thompson \cite[Theorem 4.14]{Lindo}, a ring $R$ is left von Neumann regular if and only if every $R$-module is trace in an $FP_1$-injective preenvelope up to isomorphism. Since $n$-von Neumann regular rings can be defined in terms of $FP_n$-injective modules \cite[Theorem 4.8]{MD4}, we obtain the following:

\begin{prop}\label{trace Neuamnn} Let \(n > 1\) be an integer.  A ring $R$ is left  $n$-von Neumann regular if and only if every $R$-module is trace in an $FP_n$-injective preenvelope up to isomorphism.
\end{prop}

\begin{proof}
The equivalence can be obtained by applying \cite[Theorem 4.1]{Lindo} to the categories $\mathcal{U}=\text{R}$-$\text{Mod}$ and $\mathcal{V}= \mathcal{FP}_n$-$Inj(R)$. It should be noted that, according to Remark \ref{suficientes proy}, every $R$-module can be injected into $\mathcal{V}$ and according to \cite[Proposition 3.5]{BEI}, the class of $FP_n$-injective modules is preenveloping for any ring $R$.
\end{proof}

Let $R$ be a commutative local ring and $\mathsf{m}$ its maximal ideal. By \cite[Lemma 5.3]{BHPST}, $R$ is a $2$-von Neumann regular ring if and only if $R/\mathsf{m} \in \mathcal{FP}_2$-$Inj(R)$. A direct consequence of Proposition \ref{trace Neuamnn} and this fact is the following corollary.

\begin{coro}Let $R$ be a commutative local ring and $\mathsf{m}$ its maximal ideal. Then  every $R$-module is trace in an $FP_2$-injective preenvelope up to isomorphism if and only if $R/\mathsf{m} \in \mathcal{FP}_2$-$Inj(R)$.
\end{coro}
\qed

Let $I$ be an ideal of $R$. We say that $I$ is  a \textit{trace ideal} of $R$ if $I=\tau_N(R)$ for some $R$-module $N$, as defined in \cite[Definition 2.1]{Lindo}. The result below follows directly from \cite[Theorem 4.1]{Lindo}, \cite[Theorem 5.2]{Lindo}, and Remark \ref{envelopin}. 

\begin{coro} Let $R$ be a ring  such that the class of $\mathcal{FP}_n$-$Proj(R)$ is closed under pure quotients. Then: 
\begin{enumerate}
\item $R$ is a left $n$-von Neumann regular ring if and only if every $R$-module is trace in its  $FP_n$-injective envelope up to isomorphism.
\item $R$ is a left self $FP_n$-injective ring if and only if every ideal (principal) $I\subseteq R$ with $R/I\in \mathcal{FP}_n$-$Proj(R)$ is trace in its $FP_n$-injective envelope.

Indeed, if $R$ is a left self $FP_n$-injective ring and $I$ is a ideal of $R$  with $R/I\in \mathcal{FP}_n$-$Proj(R)$, then $I$ is a trace ideal of $R$.
\end{enumerate}
\end{coro}
\qed

\section*{Acknowledgements}
The authors would like to thank anonymous referee for corrections and helpful suggestions that essentially improve the original version of the paper.

The authors gratefully  acknowledge financial support from CSIC (\textit{Comisi\'on Sectorial de Investigaci\'on Cient\'ifica}) of Uruguay (grant no. 22520220100067UD). Additionally, the research of the second author was partially supported by the \textit{Programa Despegue Cient\'ifico} of Pedeciba Matem\'atica (Uruguay) .

\end{document}